\documentclass[a4paper,12pt]{article}

\usepackage[margin=2.8cm]{geometry} 
\usepackage{fleqn} 

\usepackage{graphicx}
\usepackage{amssymb}
\usepackage{amsmath}
\usepackage{amsfonts}
\usepackage{amsthm}

\usepackage{subfigure}
\usepackage{natbib}
\usepackage{epstopdf}
\usepackage{appendix}

\newcommand{\R}{\mathbb{R}}
\newcommand{\N}{\mathbb{N}}
\newcommand{\zM}{z_{max}}
\newcommand{\zm}{z_{min}}
\newcommand{\zf}{z_{F}}
\newcommand{\zw}{z_{W}}
\newcommand{\normx}[1]{\left\|{#1}\right\|_X}

\newtheorem{assumption}{Assumption}
\newtheorem{prop}{Proposition}
\newtheorem{theorem}{Theorem}[section]
\newtheorem{definition}{Definition}
\newtheorem*{acknowledgements}{Acknowledgements}

\title{Immuno-epidemiology of a population
structured by immune status: a mathematical study of\\ waning immunity
and immune system boosting}
\author{M.~V. Barbarossa \thanks{Bolyai Institute,
University of Szeged, H-6720 Szeged, Aradi v\'{e}rtan\'{u}k tere 1, Hungary,
\textit{barbaros@math.u-szeged.hu}} \, and 
      G.~R\"ost \thanks{Bolyai Institute,
University of Szeged, H-6720 Szeged, Aradi v\'{e}rtan\'{u}k tere 1, Hungary,
\textit{rost@math.u-szeged.hu}}}

\begin{document}

\maketitle

\begin{abstract}
When the body gets infected by a pathogen the immune system develops pathogen-specific immunity. Induced immunity decays in time and years after recovery the host might become susceptible again. Exposure to the pathogen in the environment boosts the immune system thus prolonging the time in which a recovered individual is immune. Such an interplay of within host processes and population dynamics poses significant challenges in rigorous mathematical modeling of immuno-epidemiology.
We propose a new framework to model SIRS dynamics, monitoring the immune status of individuals and including both waning immunity and immune system boosting. Our model is formulated as a system of two ODEs coupled with a PDE. After showing existence and uniqueness of a classical solution, we investigate the local and the global asymptotic stability of the unique disease-free stationary solution. Under particular assumptions on the general model, we can recover known examples such as large systems of ODEs for SIRWS dynamics, as well as SIRS with constant delay. Moreover, a new class of SIS models with delay can be obtained in this framework.\\
\ \\
\textit{KEYWORDS:} {Immuno-epidemiology;\, Waning immunity;\, Immune status;\, Boosting;\, Physiological structure;\, Reinfection;\, Delay equations;\, Global stability;\, Abstract Cauchy problem}\\
\ \\
\textit{AMS Classification:} {92D30;\, 35Q91;\, 47J35;\, 34K17}
\end{abstract}

\section{Introduction}
\label{sec:intro}
Models of SIRS type are a traditional topic in mathematical epidemiology. Classical approaches present a population divided into susceptibles (S), infectives (I) and recovered (R), and consider interactions and transitions among these compartments. Susceptibles are those hosts who either have not contracted the disease in the past or have lost immunity against the disease-causing pathogen. When a susceptible host gets in contact with an infective one, the pathogen can be transmitted from the infective to the susceptible and with a certain probability, the susceptible host becomes infective himself. After pathogen clearance the infective host recovers and becomes immune for some time, afterwards he possibly becomes susceptible again.\\
\ \\
From the in-host point of view, immunity to a pathogen is the result of either active or passive immunization. The latter one is a transient protection and is due to the transmission of antibodies from the mother to the fetus through the placenta, thanks to which the newborn is immune for several months after birth \cite{McLean1988a}. Active immunization is either induced by natural infection or can be achieved by vaccine administration \cite{Siegrist2008,KubyImmBook}.

\begin{figure}[!]
\centering
\includegraphics[width=0.9\columnwidth]{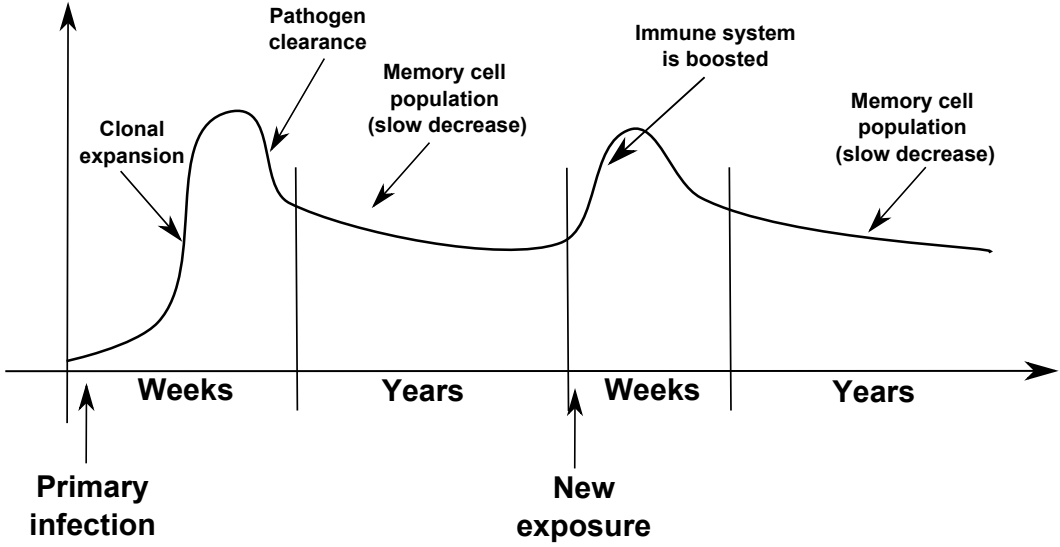}
\caption{Level of pathogen-specific immune cells with respect to the time. Generation of memory cells takes a few weeks: Once primary infection occurred the adaptive immune system produces a high number of specific immune cells (clonal expansion). After pathogen clearance, specific immune cells (memory cells) are maintained for years at a level that is much higher than in a naive host. Memory cells are activated in case of secondary infection.}
\label{Fig:introfig1}
\end{figure}

\indent Let us first consider the case of natural infection. A susceptible host, also called \textit{naive host}, has a very low level of specific immune cells for a pathogen (mostly a virus or a bacterium, but possibly also a fungus). 
The first response to a pathogen is nonspecific, as the innate immune system cannot recognize the physical structure of the pathogen. The innate immune response slows down the initial growth of the pathogen, while the adaptive (pathogen-specific) immune response is activated. Clonal expansion of specific immune cells (mostly antibodies or CTL cells) and pathogen clearance follow. The population of pathogen-specific immune cells is maintained for long time at a level that is much higher than in a naive host. These are the so-called \textit{memory cells} and are activated in case of secondary infection (see Fig. \ref{Fig:introfig1}). Memory cells rapidly activate the immune response and the host mostly shows mild or no symptoms \cite{Wodarz2007,Antia2005}.

\indent Each exposure to the pathogen might have a boosting effect on the population of specific memory cells. Indeed, the immune system reacts to a new exposure as it did during primary infection, thus yielding an increased level of memory cells \cite{Antia2005}. Though persisting for long time after pathogen clearance, the memory cell population slowly decays, and in the long run the host might lose his pathogen-specific immunity. \citet{Wodarz2007} writes that it ``is unclear for how long hosts are protected against reinfection, and this may vary from case to case.'' 

\indent  Vaccine-induced immunity works similarly to immunity induced by the natural infection. Agents contained in vaccines resemble, in a weaker form, the disease-causing pathogen and force a specific immune reaction without leading to the disease. If the vaccine is successful, the host is immunized for some time. Vaccinees experience immune system boosting and waning immunity, just as hosts recovered from natural infection do. In general, however, disease-induced immunity induces a much longer lasting protection than vaccine-induced immunity does \cite{Siegrist2008}. Waning immunity might be one of the factors which cause, also in highly developed regions, recurrent outbreaks of infectious diseases such as measles, chickenpox and pertussis. On the other side, immune system boosting due to contact with infectives prolongs the protection duration. 

\indent In order to understand the role played by waning immunity and immune system boosting in epidemic outbreaks, in the recent past several mathematical models were proposed. Few of these models describe only in-host processes which occur during and after the infection \cite{Wodarz2007,Heffernan2008,Diekmann2014}. Many more models, formulated in terms of ordinary differential equations (ODEs), consider the problem only at population level, defining compartments for individuals with different levels of immunity and introducing transitions between these compartments \cite{Dafilis2012,Heffernan2009}. Possibly, also vaccinated hosts or newborns with passive immunity are included, and waning of vaccine-induced or passive immunity are observed \cite{Rouderfer1994,Mossong1999,Moghadas2008}. In some cases, e.~g. in the model by \citet{Heffernan2009}, a discretization of the immune status leads to a very large system of ODEs. Simplified versions were suggested by \citet{Glass2003b,Grenfell2012} who add the class W of hosts with waning immunity. W-hosts can receive immune system boosting due to contact with infectives or move back to the susceptible compartment due to immunity loss. \citet{Lavine2011} extend the SIRWS model in \cite{Glass2003b}, dividing each population into
age-classes of length of 6 months each and further classifying immune hosts (R) and hosts with waning immunity (W) by the level of immunity. 

\indent To describe the sole waning immunity process, authors have chosen delay differential equations (DDEs) models with constant or distributed delay \cite{Kyrychko2005,Taylor2009,Blyuss2010,Bhat2012,Belair2013}. The delay represents the average duration of the disease-induced immunity. However, neither a constant nor a distributed delay allows for the description of immune system boosting.

\indent \citet{Martcheva2006} suggest an SIRS model in which infective and recovered hosts are structured by their immune status. In infective hosts the immune status increases over the course of infection, while in recovered hosts the immune status decays at some nonconstant rate. When the immune status has reached a critical level, recovered hosts transit from the immune to the susceptible compartment.\\
\ \\
The goal of the present paper is to suggest a general framework for modeling waning immunity and immune system boosting in the context of population dynamics. We want to understand how the population dynamics, and in particular the number of infectives, affects waning immunity and immune system boosting in a recovered host and how, in turn, these two in-host processes influence the population dynamics. In contrast to the models proposed in \cite{Heffernan2009,Lavine2011}, we shall maintain the number of equations as low as possible. For the sake of simplicity we do not include vaccination in this model. \\
\indent We suggest a model in which the immune population is structured by the level of immunity, $z \in [\zm,\zM]$. Individuals who recover at time $t$ enter the immune compartment (R) with maximal level of immunity $\zM$. The level of immunity tends to decay in time and when it reaches the minimal value $\zm$, the host becomes susceptible again.
Immune system boosting is included by the mean of the probability $p(z,\tilde z),\, z\geq \tilde z,\,z,\tilde z \in \R$, that an individual with immunity level $\tilde z$ moves to immunity level $z$, when exposed to the pathogen. At the same time we choose the susceptible and the infective populations to be non-structured. In this way we combine a PDE for the immune population with two ODEs for susceptible and infective hosts. The physiological structure in the immune population is the key innovation of the modeling approach, as it allows to describe at once both waning immunity (as a natural process which occurs as time elapses) and immune system boosting. There is hence no need to include any compartment for hosts with waning immunity as it was done, e.~g. in \cite{Glass2003b,Grenfell2012}.\\
\ \\
The paper is organized as follows. In Sect.~\ref{sec:model} we carefully derive the model equations for the SIRS dynamics. This system will be referred to as model (M1). ODEs for susceptibles and infectives are easily introduced, whereas the PDE for the structured immune population is derived from balance laws. Though at a first glance the resulting immune hosts equation might resemble a size-structured model, there is a crucial difference. In size-structured models, indeed, transitions occur only in one direction, i.~e., individuals grow in size but never shrink. On the other hand, our physiologically structured population is governed by a differential equation which includes transport (decay of immune status) and jumps (boosting to any higher immune status), describing movements to opposite directions. This makes the model analysis particularly challenging.\\
\indent In Sect.~\ref{sec:exiuni} we consider fundamental properties of solutions for the model (M1). Writing the system as an abstract Cauchy problem in the state-space 
$X:=\R\times \R \times L^1\left([\zm,\zM];\R\right) \setminus \left\{0,0,0\right\}$, 
we can guarantee existence of a unique classical solution. The proof requires a quite long computation which is postponed to the appendix. Further, using an appropriate comparison system, we show nonnegativity of solutions and determine an invariant subset in $X$.\\
\indent In Sect.~\ref{sec:dfe} we consider stationary solutions. The proof of existence of stationary solution requires connection with the theory of Volterra integral equations.
We show that there exists a unique disease free equilibrium (DFE) and investigate its local and global stability. Investigation of endemic equilibria is considered in a following work.\\

\noindent Finally in Sect.~\ref{sec:connODEs} and Sect.~\ref{sec:connDDEs} we show how to obtain from the general model (M1) various systems of ODEs, such as those in \cite{Glass2003b,Grenfell2012,Heffernan2009,Lavine2011}, and systems of equations with constant delay, such as those in \cite{Taylor2009,Kyrychko2005}. Moreover we obtain as a special case a new class of SIS systems with delay.

\section{A new class of models}
\label{sec:model}
In this section we derive the model equations for the SIRS dynamics. We start with the ODEs for susceptibles and infectives, which can be easily introduced, and continue with a PDE for the structured immune population, which we obtain from a discrete approach. The result is the general model (M1). We conclude the session with remarks on the total immune population.

\subsection{Susceptibles and infectives}
\label{sec:model_SI}
Let $S(t)$ and $I(t)$ denote the total population of susceptibles, respectively infectives, at time $t$. In contrast to previous models which include short-term passive immunity \cite{Rouderfer1994,mclean1988,Moghadas2008}, we assume that newborns are all susceptible to the disease. To make the model closer to real world phenomena we suppose that newborns enter the susceptible population at rate $b(N)$, dependent on the total population size $N$.
In general one could choose the natural death rate $d(N)$ to be a function of the total population, too. However, for simplicity of computation, we consider the case of constant death rate $d(N)\equiv d>0$. The death rate $d$ is assumed to be the same for all individuals, but infectives might die due to the infection, too.
\begin{assumption}
\label{ass:bNdN} 
Let $b:[0,\infty)\to [0, b_+],\, N\mapsto b(N),$ with $0< b_+ <\infty$, be a nonnegative $C^1$-function, with $b(0)=0$.
\end{assumption}
\noindent To guarantee that there exists a nontrivial equilibrium $N^*> 0$, such that $b(N^*)=d N^*$ (see Fig. \ref{Fig:Nstar}), we require the following.
\begin{assumption}
\label{ass:bN}
There is an $N^*>0$ such that for $N\in (0,N^*)$ we have that \linebreak $b(N)>dN$, whereas for ${N>N^*}$ we have $dN>b(N)$.
Further it holds that $${b'(N^*):=\frac{d\,b(N)}{dN}|_{N=N^*}<d}.$$
\end{assumption}
When we include disease-induced death at rate $d_I>0$, an equilibrium satisfies
\begin{equation}
\label{eq:equil_Nstart_dI}
b(N)=dN +d_I I.
\end{equation}

\noindent Contact with infectives (at rate $\beta I/N$) induces susceptible hosts to become infective themselves.  Infected hosts recover at rate $\gamma>0$, that is, $1/ \gamma$ is the average infection duration. Once recovered from the infection, individuals become immune, however there is no guarantee for life-long protection. Immune hosts who experience immunity loss become susceptible again.\\
\ \\
The dynamics of $S$ and $I$ is described by the following equations:
\begin{equation*}
\begin{aligned}\label{sys:SI_decription}
S'(t) & = \underbrace{b(N(t))}_{\mbox{birth}} -\underbrace{\beta \frac{S(t)I(t)}{N(t)}}_{\mbox{infection}}-\underbrace{dS(t)}_{\mbox{death}}+\underbrace{\Lambda}_{\mbox{immunity loss}},\\[0.5em]
I'(t) & = \underbrace{\beta \frac{S(t)I(t)}{N(t)}}_{\mbox{infection}}-\underbrace{\gamma I(t)}_{\mbox{recovery}}
-\underbrace{d I(t)}_{\substack{\text{natural}\\ \text{death}}}-\underbrace{d_I I(t)}_{\substack{\text{disease-induced}\\ \text{death}}}.
\end{aligned}
\end{equation*}
The term $\Lambda$, which represents transitions from the immune compartment to the susceptible one, will be specified below together with the dynamics of the immune population.

\begin{figure}[!]
\centering
\includegraphics[width=0.6\columnwidth]{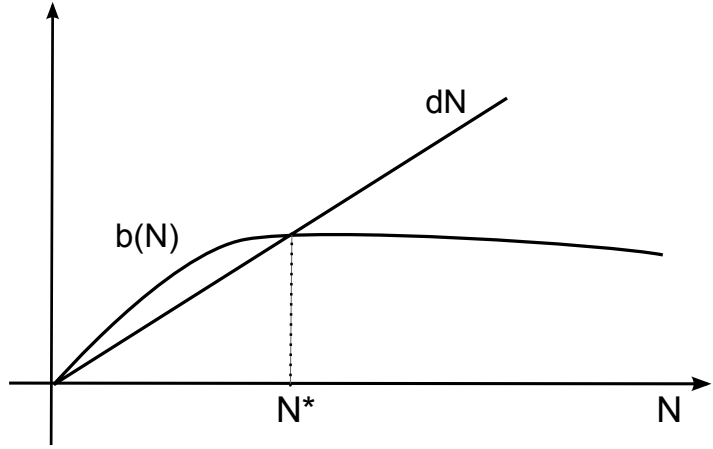}
\caption{Model Assumptions: Birth and death rate as functions of the total population $N$. In absence of disease-induced death there exists an equilibrium $N^*$ such that $b(N^*)=dN^*$.}
\label{Fig:Nstar}
\label{Fig:Nstar}
\end{figure}

\subsection{Immune hosts}
\label{sec:model_rstr}

Let $r(t,z)$ denote the density of immune individuals at time $t$ with immunity level $z \in [\zm, \zM],\,0\leq \zm <\zM<\infty$. The total population of immune hosts is given by 
\begin{equation*}
R(t) = \int_{\zm}^{\zM} r(t,z)\, dz.
\end{equation*}
The parameter $z$ describes the immune status and can be related to the number of specific immune cells of the host (cf. Sect.~\ref{sec:intro}). The value $\zM$ corresponds to maximal immunity, whereas $\zm$ corresponds to low level of immunity.

\begin{figure}[!] 
\begin{center}
\subfigure{\includegraphics[width=0.43\linewidth]{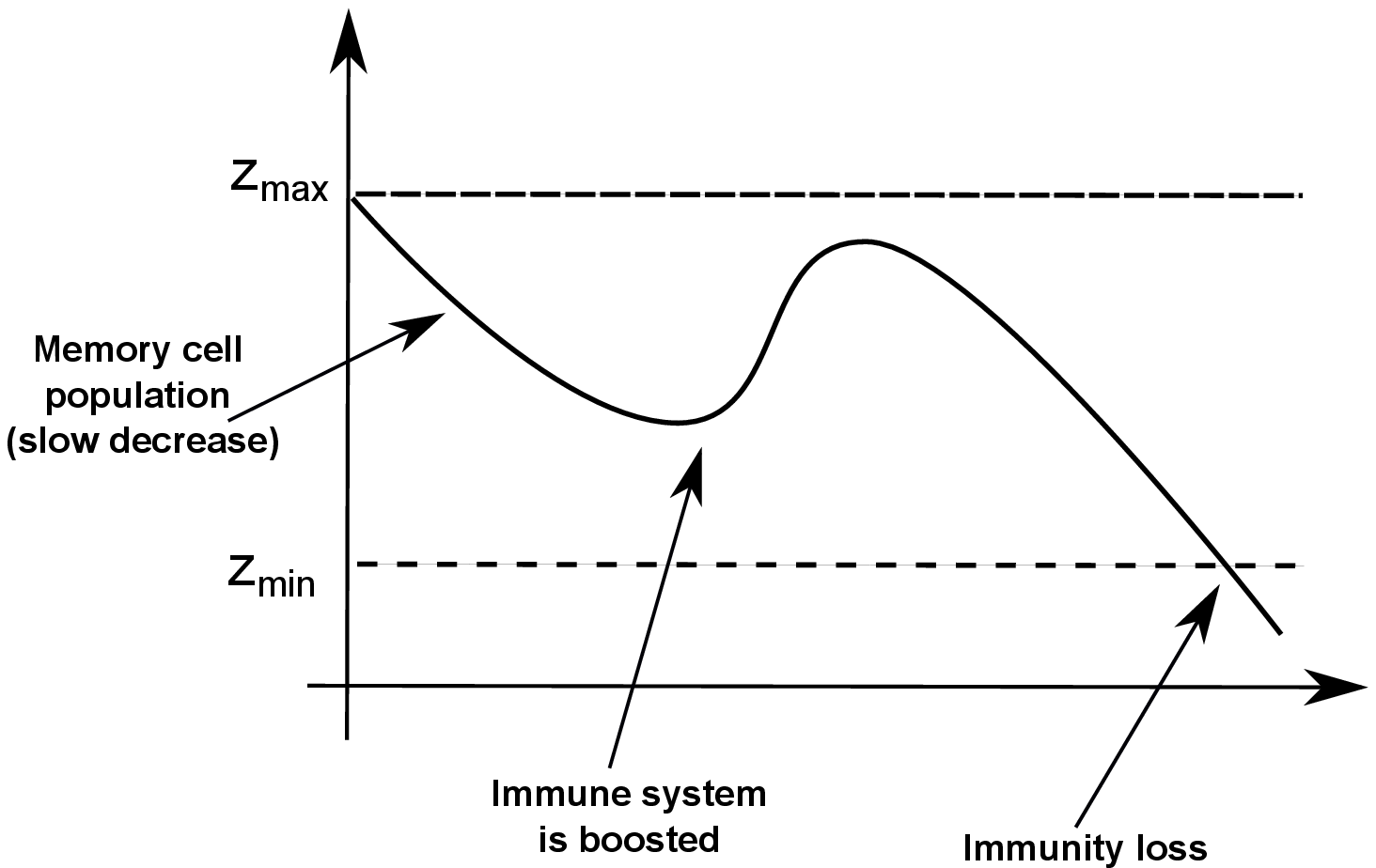}}
\subfigure{\includegraphics[width=0.55\linewidth]{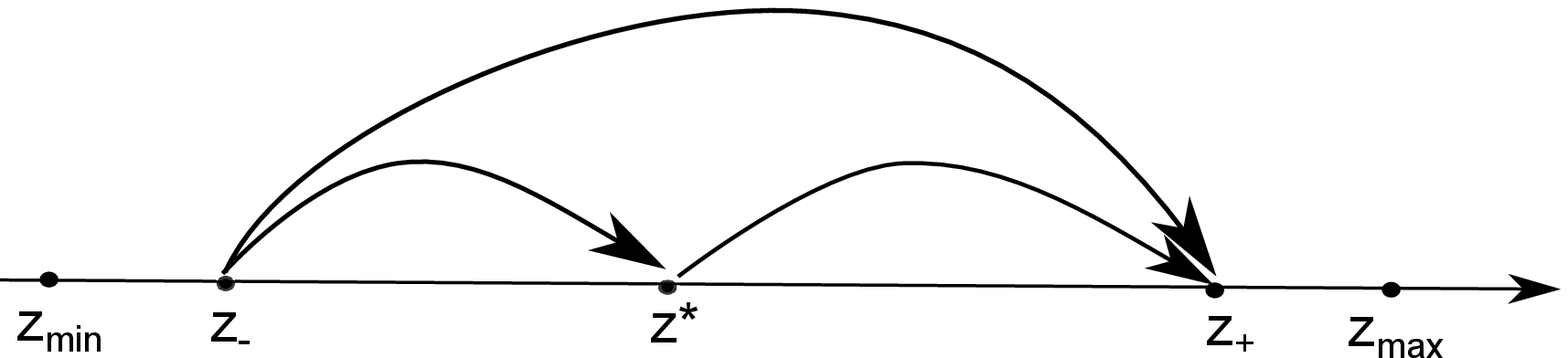}}
\caption{Model assumptions for the immune population. (a) Natural infection induces the maximal level of immunity $\zM$. The level of immunity decays in time and when it reaches the minimal value $\zm$, the host becomes susceptible again. Exposure to the pathogen has a boosting effect on the immune system. 
(b) Contact with infectives can boost the immunity level $z\in [\zm,\zM]$ to any higher value.}
 \end{center}
 \label{Fig:zminzmax}
 \end{figure}

\noindent We assume that individuals who recover at time $t$ enter the immune compartment (R) with maximal level of immunity $\zM$. The level of immunity tends to decay in time and when it reaches the minimal value $\zm$, the host becomes susceptible again. However, contact with infectives, or equivalently, exposure to the pathogen, can boost the immune system from $z\in [\zm,\zM]$ to any higher immune status, see Fig. \ref{Fig:zminzmax}. It is not straightforward to determine how this kind of immune system boosting works, as no experimental data are available. Nevertheless, laboratory analysis on vaccines tested on animals or humans suggest that the boosting efficacy might depend on several factors, among which the current immune status of the recovered host and the amount of pathogen he receives \cite{Amanna2007,Luo2012}. Possibly, exposure to the pathogen can restore the maximal level of immunity, just as natural infection does (we shall consider this special case in Sect.~\ref{sec:ex_maxlevel_boostingmodel}).\\
\indent Let $p(z,\tilde z),\, z\geq \tilde z,\,z,\tilde z \in \R$ denote the probability that an individual with immunity level $\tilde z$ moves to immunity level $z$, when exposed to the pathogen. Due to the definition of $p(z,\tilde z)$ we have 
$p(z,\tilde z)\in [0,1],\, z\geq \tilde z$ and $p(z,\tilde z)= 0$, for all $z < \tilde z$.
As we effectively consider only immunity levels in the interval $[\zm,\zM]$, we set
\begin{displaymath}
p(z,\tilde z)= 0, \quad \mbox{for all} \quad \tilde z \in (-\infty,\zm) \cup (\zM,\infty).	
\end{displaymath}
Then we have
$$\int_{-\infty}^{\infty}p(z,\tilde z)\, dz\,=\,\int_{\tilde z}^{\zM}p(z,\tilde z)\, dz\,=\,1,\quad  \mbox{for all}\quad \tilde z \in [\zm,\zM].$$
Exposure to the pathogen might restore exactly the immunity level induced by the disease ($\zM$). In order to capture this particular aspect of immune system boosting, we write the probability $p(z,\tilde z)$ as the combination of a continuous ($p_0$) and atomic measures (Dirac delta):
\begin{displaymath}
p(z,\tilde z)= c_{max}(\tilde z)\delta (\zM-\tilde z) + c_0(\tilde z)p_0(z,\tilde z) + c_1(\tilde z)\delta(z-\tilde z),
\end{displaymath}
where 
\begin{description}
  \item[{$c_{max}:[\zm,\zM]\to [0,1],\;y\mapsto c_{max}(y)$},] is a continuously differentiable function and describes the probability that, due to contact with infectives, a host with immunity level $y$ boosts to the maximal level of immunity $\zM$.	
	\item[{$c_{0}:[\zm,\zM]\to [0,1],\;y\mapsto c_{0}(y)$},] is a continuously differentiable function and describes the probability that, due to contact with infectives, a host with immunity level $y$ boosts to any other level $z \in (y,\zM)$, according to the continuous probability $p_0(z,y)$. 
	\item[{$c_{1}(y):[\zm,\zM]\to [0,1],\;y\mapsto c_{1}(y)=1-c_{max}(y)-c_0(y)$},] describes the	
	probability that getting in contact with infectives, the host with immunity level $y$ does not experience immune system boosting. 
	\end{description}

\noindent The immunity level decays in time at some rate $g(z)$ which is the same for all immune individuals with immunity level $z$, that is,
\begin{equation*}
\frac{d}{dt}z(t)=g(z).
\end{equation*}

\begin{assumption}
\label{ass:g}
Let $g:[\zm,\zM]\to (0,K_g],\; K_g<\infty$ be continuously differentiable.
\end{assumption}
The positivity of $g(z)$ is given by the fact that, if  $g(\tilde z)=0$ for some value $\tilde z \in [\zm,\zM]$, there would be no change of the immunity level at $\tilde z$, contradicting the hypothesis of natural decay of immune status.\\
\ \\
In absence of immune system boosting, an infected host who recovered at time $t_0$ becomes again susceptible at time $t_0+T$, where
$$ T=\int_{\zm}^{\zM} \frac{1}{g(z)}\,dz.$$
The above equality becomes evident with an appropriate change of variables. Setting $u=z(s)$, we find 
$du=z'(s)\,ds=g(z(s))\,ds$. For the integration boundaries we have $\zM=z(t_0)$ and $\zm=z(t_0+T)$. It follows that
$$ \int_{\zm}^{\zM} \frac{1}{g(z)}\,dz= \int_{t_0}^{t_0+T}ds=T.$$

\noindent To obtain a correct physical formulation of the equation for $r(t,z)$, we start from a discrete approach, as it could be done for age-structured or size-structured models \cite{WebbPopChap,Ellner2009}. The number of immune individuals in the immunity range $[ z -\Delta z, z],\, z \in [\zm,\zM],$ is $r(t,z)\Delta z$.\\
\indent We track how many individuals enter and exit a small immunity interval $[z -\Delta z,z]$ in a short time $\Delta t\ll 1$. 
Starting at a time $t$ with $r(t,z)\Delta z$ individuals, we want to describe the number of individuals at time $t+\Delta t$. 
Recall that the immunity level tends to decay, hence individuals enter the interval from $z$ and exit from $z -\Delta z$ (i.~e., transition occurs the other way around than in age-structured or size-structured models). As contact with infective individuals might boost the immune system, an immune host with immunity level in $[z -\Delta z, z]$ can move from this interval to any higher immunity level. For the same reason, immune individuals with any lower immunity level can ``jump'' to the interval $[z -\Delta z,z]$. We assume that contacts between infectives and immune hosts occur at the same rate, $\beta I/N$, as between infectives and susceptibles.\\
\indent Given any $z\in [\zm,\zM]$, denote by $G_z$ the partition of the interval $[\zm,z]$ into small intervals of length $\Delta z$. Observe $G_z \subset G_{\zM}$ for all $z\in [\zm,\zM]$.\pagebreak

\noindent The total number of immune individuals with immunity level in $[z -\Delta z, z]$ at time $t+\Delta t$ is given by
\begin{align*}
r(t+\Delta t,z)\Delta z & = \underbrace{r(t,z)\Delta z}_{\mbox{at time } t} + \underbrace{r(t,z)g(z) \Delta t}_{\mbox{incoming (waning)}}
 -\underbrace{r(t, z -\Delta z)g(z -\Delta z)\Delta t}_{\mbox{outgoing (waning)}}\\[0.5em]
& \phantom{=}- \underbrace{dr(t,z)\Delta z\Delta t}_{\mbox{die out}} -\underbrace{\beta \frac{I(t)}{N(t)} r(t,z)\Delta z \Delta t}_{\mbox{outgoing (boosting)}} \\[0.5em]
  & \phantom{=} +\underbrace{\sum_{v \in G_z}\beta \frac{I(t)}{N(t)} p(z,v)\Delta z\,r(t,v)\Delta v\,\Delta t}_{\mbox{incoming (boosting)}}.												
\end{align*}
Perform division by $\Delta z>0$, 
\begin{align*}
r(t+\Delta t,z) & = r(t,z) + r(t,z)g(z) \frac{\Delta t}{\Delta z} -r(t,z -\Delta z)g(z -\Delta z)\frac{\Delta t}{\Delta z}- dr(t,z)\Delta t\\[0.6em]
  & \phantom{=} -\beta \frac{I(t)}{N(t)} r(t,z)\Delta t	+\sum_{v \in G_z}\beta \frac{I(t)}{N(t)} p(z,v)\,r(t,v)\Delta v\,\Delta t,											
\end{align*}
and compute the limit $\Delta z \to 0$ (observe that also $\Delta v \to 0$, as $\Delta z$ and $\Delta  v$ are elements of the same partition) to get
\begin{align*}
r(t+\Delta t,z)-r(t,z) & =  \Delta t\frac{\partial}{\partial z}\bigl(g(z)r(t,z)\bigr)- dr(t,z)\Delta t\\[0.5em]
												& \phantom{=} -\beta \frac{I(t)}{N(t)} r(t,z)\Delta t+\beta \frac{I(t)}{N(t)}\Delta t\int_{\zm}^{z}p(z,v)r(t,v)\, dv.										
\end{align*}
Finally, we divide by $\Delta t$ and let $\Delta t \to 0$. From the discrete approach derivation it becomes clear that the quantity $\Lambda$, initially introduced in the $S$ equation at p. \pageref{sys:SI_decription} to represent the hosts who experience immunity loss, is given by the number $g(\zm)r(t,\zm)$ of immune hosts who reach the minimal level of immunity.\\
\ \\
Altogether for $t\geq 0$ we have the system of equations
\begin{equation}
\begin{aligned}
S'(t) & = b(N(t)) -\beta \frac{S(t)I(t)}{N(t)}-dS(t)+\underbrace{g(\zm)r(t,\zm)}_{=:\Lambda},\\[0.5em]
I'(t) & = \beta \frac{S(t)I(t)}{N(t)}-(\gamma+d+d_I) I(t),
\end{aligned}
\label{sys:mod1_PDE}
\end{equation}
with initial values $S(0)=S^0>0,\,I(0)=I^0\geq 0$, coupled with a partial differential equation for the immune population,
\begin{equation}
\frac{\partial }{\partial t}r(t,z) - \frac{\partial }{\partial z}\left(g(z)r(t,z)\right) = -dr(t,z)
+ \beta \frac{I(t)}{N(t)}\left(\int_{\zm}^{z} p(z,v)r(t,v)\,dv- r(t,z)\right),
\label{eq:mod1_pde_R} 
\end{equation}
where $z \in [\zm,\zM]$, with boundary condition
\begin{equation}
\label{eq:BC_mod1_pde_R} 
g(\zM)r(t,\zM)= \gamma I(t) + \beta \frac{I(t)}{N(t)}\int_{\zm}^{\zM} p(\zM,v)r(t,v)\,dv,
\end{equation}
and initial distribution $r(0,z)=\psi(z),\, z \in  [\zm,\zM]$.\\
\ \\
At a first glance equation \eqref{eq:mod1_pde_R} might resemble a size-structured model, but there is an important difference. In size-structured models transitions occur only in one direction, i.~e., individuals grow in size but never shrink. The PDE \eqref{eq:mod1_pde_R}  for the immune population is governed by a transport process (decay of immune status) and jumps (boosting to any higher immune status). The model analysis results much more challenging than the one of classical size-structured models.\\
\ \\
In the following we refer to the complete system \eqref{sys:mod1_PDE} -- \eqref{eq:BC_mod1_pde_R} as to \textbf{model (M1)}. An overview on the model notation is provided in Table \ref{Tab:not_mod1pde}. 

\begin{table}
\begin{tabular}{ c | l}
\textbf{Symbol} & \textbf{Description}\\
\hline
$S(t)$ & number of susceptibles at time $t$  \\
$I(t)$ & number of infective hosts at time $t$  \\
$R(t)$ & number of immune individuals at time $t$  \\
$r(t,z)$ & density of immune individuals with immunity level $z$ at time $t$  \\
$N(t)$& total population ($N(t)=S(t)+I(t)+R(t)$) at time $t$\\
$b(N)$ & recruitment rate\\
$d$ & natural death rate\\
$d_I$ & disease-induced death rate\\
$\beta$ & transmission rate\\
$\gamma$  & recovery rate\\
$g(z)$  & natural decay of immunity\\
$p(z,\tilde z)$ & probability that boosting raises immunity level $\tilde z$ to level $z$\\
$\zM$ & maximal level of immune status in immune hosts\\
$\zm$ & minimal level of immune status in immune hosts\\
$\psi(z)$ & initial distribution for $r(t,z)$\\
\hline
\end{tabular}
\caption{Notation for model (M1). All quantities are nonnegative.}
\label{Tab:not_mod1pde}
\end{table}

\subsection{Total number of immune individuals}
\label{sec:relR}
From the PDE \eqref{eq:mod1_pde_R} -- \eqref{eq:BC_mod1_pde_R} we obtain the total number $R(t)$ of immune individuals at time $t\geq0$.
Integrating the left-hand side of \eqref{eq:mod1_pde_R} in $[\zm,\zM]$, we find
\begin{equation*}
\begin{aligned}
& \int_{\zm}^{\zM}\frac{\partial }{\partial t}r(t,z)\,dz-\int_{\zm}^{\zM}\frac{\partial }{\partial z}\left(g(z)r(t,z)\right)\,dz\\
& = \frac{\partial }{\partial t}\int_{\zm}^{\zM}r(t,z)\,dz-g(\zM)r(t,\zM)+g(\zm)r(t,\zm)\\
& = R'(t) - \underbrace{\left(\gamma I(t) + \beta \frac{I(t)}{N(t)}\int_{\zm}^{\zM} p(\zM,v)r(t,v)\,dv\right)}_{\underset{\eqref{eq:BC_mod1_pde_R}}{=} \, g(\zM)r(t,\zM)} +g(\zm)r(t,\zm).
\end{aligned}
\end{equation*}
Similarly, integrating the right-hand side of \eqref{eq:mod1_pde_R} we have
\begin{equation*}
\begin{aligned}
& -\int_{\zm}^{\zM}dr(t,z)\,dz+ \beta \frac{I(t)}{N(t)}\left[\int_{\zm}^{\zM}\left(\int_{\zm}^{z} p(z,v)r(t,v)\,dv - r(t,z)\right)\,dz\right]\\
& = -dR(t)+ \beta \frac{I(t)}{N(t)}\biggl[\int_{\zm}^{\zM}\int_{\zm}^{z} p(z,v)r(t,v)\,dv\,dz - R(t)\biggr].
\end{aligned}
\end{equation*}
Altogether:
\begin{equation}
\begin{aligned}
R'(t) & = \gamma I(t)-dR(t) -g(\zm)r(t,\zm)+\beta \frac{I(t)}{N(t)}\biggl[\int_{\zm}^{\zM} p(\zM,v)r(t,v)\,dv\\
& \phantom{==} + \int_{\zm}^{\zM}\int_{\zm}^{z} p(z,v)r(t,v)\,dv\,dz- R(t)\biggr].
\label{dotRt_rel_2a}
\end{aligned}
\end{equation}
When an immune host comes in contact with infectives, his immune system gets boosted so that either the maximal level of immunity or any other higher (or equal) level of immunity is restored. This means that the terms in the square bracket in \eqref{dotRt_rel_2a} cancel out and the total immune population at time $t$ satisfies
\begin{equation}
R'(t)= \gamma I(t)  -g(\zm)r(t,\zm) -dR(t).
\label{dotRt_rel_1}
\end{equation}
In other words, inflow at time $t$ into the immune class occurs by recovery of infected hosts, while outflow is either due to death or to immunity loss. Observe that the result \eqref{dotRt_rel_1} holds also in the special case of no boosting ($c_0(z)\equiv 0$ and $c_{max}(z)\equiv 0,\, \in [\zm,\zM])$, and in the case in which boosting always restores the maximal level of immunity ($c_0(z)\equiv 0$ and $c_{max}(z)\equiv 1,\,z \in [\zm,\zM]$).


\section{Existence, uniqueness and positivity}
\label{sec:exiuni}
Consider model (M1), with given initial data $S(0)=S^0\geq 0,\;I(0)=I^0\geq 0$ and $\psi(z)\geq 0$ for all $z \in  [\zm,\zM]$.
Let $Y:=L^1\left([\zm,\zM];\R\right)$ with its usual norm, and let $\tilde r \in Y$ be defined by
\begin{displaymath}
\tilde r(t):[\zm,\zM]\to \R,\; z\mapsto r(t,z), \qquad t\geq 0.
\end{displaymath}
For given values of $I$ and $N$, we define the map $\mathcal{F}^{(I,N)}:Y \to Y,\; \phi\mapsto \mathcal{F}^{(I,N)}\phi$ by
\begin{equation}
\label{eq:abst_cauchyprobl_F}
\left(\mathcal{F}^{(I,N)}\phi\right)(z) = \beta\frac{I}{N}\left[\int_{\zm}^{z} \phi(v)(\pi(v)(z))\,dv-\phi(z)\right],
\end{equation}
with $\pi(v)(z)=p(z,v) \in [0,1],\;v\in [\zm,\zM]$. In this notation the initial condition for $r$ is given by $\tilde r(0)(z)=\psi(z),\;z\in[\zm,\zM]$.\\
\ \\
The state space for model (M1) is $X:=\R\times \R \times Y \setminus \left\{0,0,0\right\}$ with the norm $\normx{\cdot}$ defined by
\begin{displaymath}
	\normx{x}:=|x_1|+|x_2|+\int_{\zm}^{\zM}|x_3(z)|\,dz,\quad \mbox{for all}\quad x=(x_1,x_2,x_3)\in X.
\end{displaymath}
Observe that for well-posedness of the model (M1) the zero is not an element of $X$. For all $t\geq 0$, we define $x(t)=(S(t),I(t),\tilde r(t)) \in X$. Then the model (M1) can be formulated as a perturbation of a linear abstract Cauchy problem,
\begin{equation}
\begin{aligned}
\frac{d}{dt}x(t) + A x(t) & = Q(x(t)), \qquad t>0,\\
x(0)&=(S^0,I^0,\psi),
\end{aligned}
\label{abst_cauch_probl_X}
\end{equation}
where $A$ and $Q$ are respectively a linear and a nonlinear operator on $X$, defined by
\begin{align*}
A_1(x_1,x_2,x_3) &: = dx_1+g(\zm)x_3(\zm),\\[0.3em]
A_2(x_1,x_2,x_3) &: = (\gamma+d+d_I)x_2,\\[0.3em]
A_3(x_1,x_2,x_3) &: = d x_3 + \frac{\partial }{\partial z}\left( g(z)x_3(z)\right),
\end{align*}
and 
\begin{align*}
Q_1(x_1,x_2,x_3) &: = b(\hat x)-\beta \frac{x_1x_2}{\hat x},\\[0.3em]
Q_2(x_1,x_2,x_3) &: = \beta \frac{x_1x_2}{\hat x},\\[0.3em]
Q_3(x_1,x_2,x_3) &: = \mathcal{F}^{(x_2,\hat x)} x_3,
\end{align*}
with $\hat x:=x_1+x_2+\int_{\zm}^{\zM}x_3(z)\,dz$.\\
\ \\
The next result guarantees the existence and uniqueness of a \textit{classical solution}. We remark that it is called \textit{mild solution} of \eqref{abst_cauch_probl_X} the continuous solution $x$ of the integral equation
\begin{equation*}
x(t)= T(t)u(0)+\int_0^t T(t-s)Q(x(s))\,ds,
\end{equation*}
where $\left\{T(t)\right\}_{t\geq0}$ is the $\mathcal C_0$-semigroup on $X$ generated by $-A$. On the other side, a function $x:[0,T)\to X$ is a \textit{classical solution} of \eqref{abst_cauch_probl_X} on $[0,T)$ if $x$ is continuous on $[0,T)$, continuously differentiable on $(0,T)$, $u(t)$ is in the domain $D(A)$ of $A$ for all $t \in (0,T)$ and \eqref{abst_cauch_probl_X} is satisfied on $[0,T)$ (see \cite{Pazy1983}). 

\begin{theorem}
\label{thm:existence_M1_solux}
Let $S(0)=S^0\geq 0,\;I(0)=I^0\geq 0$ and $\psi(z)\geq 0$ for all $z \in  [\zm,\zM]$ be given. Let $d>0,\;b:[0,\infty)\to [0, b_+],\,0< b_+ <\infty$ satisfy Assumption \ref{ass:bNdN} and Assumption \ref{ass:bN}, and $g:[\zm,\zM]\to (0,K_g),\; 0<K_g<\infty$ satisfy Assumption \ref{ass:g}. For given values $I$ and $N$, let $\mathcal F^{(I,N)}$ be defined by \eqref{eq:abst_cauchyprobl_F}.\\
\indent Then there exists a unique solution $x$ on $[0,\infty)$ of the abstract Cauchy problem \eqref{abst_cauch_probl_X}, with initial data $x(0)=(S^0,I^0,\psi) \in X$.
\end{theorem}
\ \\
\textit{Proof.} To have existence of a unique (classical) solution on $[0,\infty)$ we need to show:
(i) that $-A$ is the generator of a $\mathcal C_0$-semigroup on $X$ and (ii) that $Q$ is continuously differentiable in $X$ (see \cite[Chap.~6]{Pazy1983}).\\
\ \\
(i) The first hypothesis is easily verified as $A$ corresponds to the linear homogeneous part of the system and its domain is 
\begin{align*}
D(A) & = \bigl\{ x \in \R \times \R \times C^1([\zm,\zM])\quad \mbox{such that}\\
     & \phantom{= \bigl\{} A_3(x) \in Y \;\mbox{and}\; g(\zM)x_3(\zM)=\gamma x_2\bigr\} \subset X.
\end{align*}
Similar linear operators arising from population dynamics were considered in \cite{WebbPopChap,CalsinaFarkas2012}.\\
\ \\
(ii) Continuous differentiability of $Q$ can be shown in two steps. First, for all $x,\,w \in X$, we determine the existence of the operator $DQ(x;w)$ defined by  
\begin{displaymath}
DQ(x;w):=\lim \limits_{h \to 0}\frac{Q(x+hw)-Q(x)}{h}.
\end{displaymath}
Second, we show that the operator $DQ(x;\cdot)$ is continuous in $x$, that is
\begin{displaymath}
\lim \|DQ(x;\cdot) -DQ(y;\cdot) \|_{OP}=0 \qquad \mbox{for} \qquad \normx{x-y}\to 0,
\end{displaymath}
where $ \| \cdot \|_{OP}$ is the operator norm. These two steps require the long computation included in the Appendix. \qed 

\noindent Note that to have the existence of a classical solution one has to show continuous differentiability of $Q$. To have existence and uniqueness of a mild solution, as well as continuous dependence on the initial data, it is sufficient to prove Lipschitz continuity of $Q$ \citep[see][]{Pazy1983}.\\
\ \\
\noindent From now on, we shall assume that all hypotheses of Theorem \ref{thm:existence_M1_solux} hold.\\
\ \\
Next we show that, given nonnegative initial data, model (M1) preserves nonnegativity. We proceed in two steps. First we introduce model (M2) in which boosting restores the maximal immune status. The PDE in model (M2) has the same characteristic curves as the PDE in model (M1), allowing us to use it as a comparison system in the proof of nonnegativity. There are other reasons to consider model (M2) separately. On the one side, the assumption of boosting restoring maximal immunity has been used in previous models, e.g. in \cite{Grenfell2012}. On the other side, in Sect.~\ref{sec:newsis} we shall use model (M2) to obtain a new SIRS system with delay.

\subsection{Model (M2): Boosting restores the maximal level of immunity}
\label{sec:ex_maxlevel_boostingmodel}
Let us assume that whenever an immune host comes in contact with the pathogen, his immune system is boosted in such a way that the maximal immunity level is restored. This means that
$$ c_{max}(z)\equiv 1,\; c_{0}(z)\equiv 0,\; z\in [\zm,\zM],$$
or equivalently,
$$ p(\zM,z)=1\; \mbox{and}\;  p(\tilde z,z)=0, \quad \mbox{for all} \quad z\in [\zm,\zM],\, \tilde z <\zM.$$
This assumption modifies the equation \eqref{eq:mod1_pde_R} and the boundary condition \eqref{eq:BC_mod1_pde_R} in model (M1) as follows
\begin{align}
\frac{\partial }{\partial t}r(t,z)-\frac{\partial }{\partial z}\left(g(z)r(t,z)\right) & = 
-dr(t,z) - r(t,z)\beta\frac{I(t)}{N(t)}, \label{sys:mod2_r1}\\
g(\zM)r(t,\zM)&= \gamma I(t) + \beta \frac{I(t)}{N(t)}R(t).
\label{sys:mod2_r2}
\end{align}
The equations for $S$ and $I$ in model (M1) remain unchanged. We shall refer to the system 
\eqref{sys:mod1_PDE}, \eqref{sys:mod2_r1} -- \eqref{sys:mod2_r2} as to \textbf{model (M2)}.\\
\ \\
Just for a moment, assume that for some $t\geq 0$ the values $I(t)$ and $S(t)$ are known, and recall $N(t)=I(t)+S(t)+R(t)$. For all $t\geq 0$ let us define
\begin{displaymath}
B(t)= \gamma I(t) + \beta \frac{I(t)}{N(t)}R(t),
\end{displaymath}
and
\begin{displaymath}
\mu(t,z)=d-g'(z)+\beta\frac{I(t)}{N(t)}.
\end{displaymath}

\begin{definition}[cf. \cite{Calsina1995}]
Let $T>0$. A nonnegative function $r(t,z)$, with $r(t,\cdot)$ integrable, is a solution of the problem \eqref{sys:mod2_r1} on $[0,T)\times[\zm,\zM]$ if the boundary condition \eqref{sys:mod2_r2} and the initial condition $r(0,z)=\psi(z),\,z\in[\zm,\zM]$ are satisfied and
\begin{displaymath}
Dr(t,z) =  -\mu(t,z)r(t,z), \qquad t\in[0,T),\,z\in \,[\zm,\zM],
\end{displaymath}
with
\begin{displaymath}
Dr(t,z): = \lim \limits_{h\to 0} \frac{r(t+h,\varphi(t+h;t,z))-r(t,z)}{h},
\end{displaymath}
where $\varphi(t;t_0,z_0)$ is the solution of the differential equation $z'(t)=g(z(t))$ with initial value $z(t_0)=z_0$.
\end{definition}
We introduce the characteristic curve $\zeta(t)=\varphi(t;0,\zM)$ which identifies the cohort of individuals who recovered (hence have maximal level of immunity) at time $t=0$. In this way we distinguish between those individuals who recovered before time $t=0$ and are already immune at the initial time, and those who recovered later than $t=0$.
Then the problem \eqref{sys:mod2_r1}--\eqref{sys:mod2_r2} can be solved along the characteristics \citep[see, e.g][]{Calsina1995,WebbPopChap} and we have the solution
\begin{equation}
\label{sol_r_M2}
r(t,z)=\begin{cases}
\frac{B(t^*)}{g(\zM)} \exp \left( -\int_{t^*}^{t}\mu(s,\varphi(s;t^*,\zM))\,ds\right),& \qquad z\geq \zeta(t)\\
& \\
\psi(\varphi(0;t,z))\exp \left( -\int_{0}^{t}\mu(s,\varphi(s;t,z))\,ds\right),& \qquad z< \zeta(t),
\end{cases}
\end{equation}
where the time $t^*$ is implicitly given by 
\begin{displaymath}
\varphi(t;t^*,\zM)=z.
\end{displaymath}
As the death rate $d>0$ is bounded, we can extend the solution to all positive times $t>0$ \citep[see also][Sec. 3.3]{Calsina1995}. It is obvious that the solution $r(t,z)$ is nonnegative for all $t\geq 0,\,z\in \,[\zm,\zM]$.

\subsection{Nonnegative solutions of model (M1)}
\label{sec:nonnegsol}
Define the set
\begin{align*}
D & =\left\{(S,I, \tilde r)\; \mbox{such that}\; S\geq 0,\; I\geq 0,\;\mbox{and}\; \tilde r(z)\geq 0,\,z\in [\zm,\zM] \right\}\; \subset X.
\end{align*}

\begin{theorem}
\label{thm:invariaceD}
The cone $D$ is positively invariant for the model (M1).
\end{theorem}
\textit{Proof.} We start with the infective population. Let $I(0)\geq 0$ be given. If $I(\bar t)=0$ for some $\bar t>0$, we have $I'(\bar t)=0$, hence the $I$ component is always nonnegative. Further we have that
\begin{equation*}
\left(\frac{I(t)}{N(t)}\right)'=\frac{I(t)}{N(t)}\left(\beta S(t)-(\gamma+d+d_I)+\frac{1}{N(t)}(b(N(t))-dN(t)-d_II(t)\right).
\end{equation*}
It follows that, given positive initial values, the total population is larger than zero for all $t>0$.\\
\ \\
The equation for $S$ includes the term $g(\zm)r(t,\zm)$, which is given by the solution $r(t,z)$ of the PDE \eqref{eq:mod1_pde_R}--\eqref{eq:BC_mod1_pde_R}. Let $S(0)\geq 0$ be given. Assume $r(t,\zm)\geq 0$ for all $t\geq 0$, hence $g(\zm)r(t,\zm)\geq 0$ (recall $g(z)>0$ by Assumption~\ref{ass:g}). We see that $S=0$ implies 
$$S'(t)= b(N(t))+g(\zm)r(t,\zm)\geq 0.$$
Hence, if $g(\zm)r(t,\zm)\geq 0$, also $S$ does not leave the nonnegative cone.\\
\ \\
To conclude the proof, we have to show that $r(t,z)\geq 0$ for all $t\geq 0$, $z\in [\zm,\zM]$. 
First we show that strictly positive initial data $r(0,z)=\psi(z)>0,\;z\in [\zm,\zM]$ lead to a strictly positive solution for all $t>0$, $z\in [\zm,\zM]$. Assuming the contrary, there is a time $\bar t>0$ such that

(i) $I(t)>0,\; S(t)\geq 0,\; N(t)>0,$ for all $t \in [0,\bar t]$;

(ii) for all $t \in [0,\bar t), \, z \in [\zm,\zM],\; r(t,z)>0$, whereas $r(\bar t, \bar z)=0$, for some $\bar z \in [\zm,\zM)$.
\\ In other words, $\bar t$ is the first time at which for some value $\bar z \in [\zm,\zM)$ the solution $r(t,z)$ of \eqref{eq:mod1_pde_R}--\eqref{eq:BC_mod1_pde_R} is zero. Note that the PDE \eqref{eq:mod1_pde_R} in model (M1) and the PDE \eqref{sys:mod2_r1} have the same characteristic curves. By assumptions (i)-(ii), along the characteristics we have the estimate
\begin{align*}
 {\frac{\partial }{\partial t}r(t,z)-\frac{\partial }{\partial z}\left(g(z)r(t,z)\right)}
& = -\left(d+ \beta\frac{I(t)}{N(t)}\right)r(t,z)\\[0.5em]
& \phantom{==} + \underbrace{\beta \frac{I(t)}{N(t)}\int_{\zm}^{z}  p(z,u)r(t,u)\,du}_{\geq 0}\\
& \geq -\left(d+ \beta\frac{I(t)}{N(t)}\right)r(t,z). 
\end{align*}
We can use the equations \eqref{sys:mod2_r1}--\eqref{sys:mod2_r2} as comparison system for \eqref{eq:mod1_pde_R}--\eqref{eq:BC_mod1_pde_R} along characteristics. From \eqref{sol_r_M2} it is evident that the solution $r(\bar t,\bar z)$ of \eqref{sys:mod2_r1}--\eqref{sys:mod2_r2} is going to be zero if and only if the characteristic curve associated to the point $(\bar t,\bar z)$ has a starting value (either $\psi(\varphi(0;\bar t,\bar z))$ or $I(t^*)$, $t^*<\bar t$) equal to zero, contradicting to (i)-(ii). We conclude that, given strictly positive initial data $\psi(z)>0$, the solution $r$ of \eqref{eq:mod1_pde_R}--\eqref{eq:BC_mod1_pde_R} is strictly positive.\\
\ \\
\noindent To complete the proof, we extend the result to nonnegative initial data. Let $\psi(z)\geq 0$. We introduce a value $\epsilon>0$ and repeat the same argument as above for initial data $\psi_\epsilon(z):=\psi(z)+\epsilon>0$. Finally we let $\epsilon\to 0$. From the continuous dependence on initial data (cf. Theorem \ref{thm:existence_M1_solux} and Pazy, 1983, Chap. 6), it follows that $r(t,z)\geq 0$ for $\psi(z)\geq 0$. \qed

\section{The disease free equilibrium}
\label{sec:dfe}
For investigation of stationary solutions of model (M1) we set the time derivative equal to zero and consider the problem
\begin{align}
0 & = b(N_*) -\beta \frac{S_*I_*}{N_*}-dS_*+g(\zm)\bar r(\zm), \label{sys:full_equi1} \\[0.3em]
0 & = \beta \frac{S_*I_*}{N_*}-(\gamma+d+d_I)I_*,\label{sys:full_equi2} \\[0.3em]
\frac{d}{d z}\left(g(z)\bar r(z)\right) & = 
d\bar r(z)- \beta \frac{I_*}{N_*}\int_{\zm}^{z}  p(z,u)\bar r(u)\,du+ \bar r(z)\beta\frac{I_*}{N_*}, \label{sys:full_equi3} \\[0.1em]
g(\zM)\bar r(\zM)&= \gamma I_* + \beta \frac{I_*}{N_*}\int_{\zm}^{\zM} p(\zM,u)\bar r(u)\,du,
\label{sys:full_equi4} 
\end{align}
where the star denotes a fixed point and the bar a stationary distribution. A stationary solution of model (M1) is a triple $\left\{S_*,I_*,\bar r(\cdot)\right\} \in X$ which satisfies \eqref{sys:full_equi1}--\eqref{sys:full_equi4}.

The total population $N_*=S_*+I_*+\int_{\zm}^{\zM} \bar r (u) du$ satisfies condition \eqref{eq:equil_Nstart_dI}.
From the equation \eqref{sys:full_equi2} we see that either $I_*=0$ or $I_*>0$ and $\beta \frac{S_*}{N_*}=(\gamma+d+d_I)$. In this paper we shall consider only the disease-free equilibrium (DFE), that is, the case $I_*=0$. The study of the endemic case $(I_*>0)$ requires a long and nontrivial analysis, including the theory of Volterra integral equations, which is beyond the scope of this manuscript.\\
\ \\
Before presenting further results we introduce the basic reproduction number $R_0$ of model (M1),
\begin{equation*}
R_0=\frac{\beta}{\gamma+d+d_I},
\end{equation*}
which indicates the average number of secondary infections generated in a fully susceptible population by one infected host over the course of his infection. 
The basic reproduction number is a reference parameter in mathematical epidemiology used to understand if, and in which proportion, the disease will spread among the population.

\begin{prop}[Existence and uniqueness of DFE]
There is exactly one disease free equilibrium (DFE), namely $(S_*,0,\bar r(\cdot))=(N^*,0,0)$.
\end{prop}
\textit{Proof.} In the case $I_*=0$, from \eqref{sys:full_equi3}--\eqref{sys:full_equi4} we have
\begin{align*}
\frac{d}{dz} \bar r(z) &= \frac{1}{g(z)}\left(d- g'(z)\right) \bar r(z) ,\\
g(\zM)\bar r(\zM) &= 0,
\end{align*}
hence the trivial solution $\bar r(z) \equiv 0$. Since $\bar r(\zm)=0$, from \eqref{sys:full_equi1} it follows that
${b(N_*)-dS_*=0}$. In particular we have $R_*=\int_{\zm}^{\zM}\bar r(z)\,dz=0$ and $S_*=N_*$. From Assumption 2 and condition \eqref{eq:equil_Nstart_dI}, we obtain $N_*=N^*$.\qed\\
\ \\

\begin{theorem}[Local stability of DFE]
If $R_0<1$, the DFE is locally stable.  \label{thmDFElocstab}
\end{theorem}
\textit{Proof.} We prove the stability from first principle. Fix $\varepsilon >0$. From the Assumptions \ref{ass:bNdN} and \ref{ass:bN} on $b(N)$ and $d$, there exists an $\omega>0$ such that $\omega<\varepsilon/3$ and 
$b(N)-dN$ is monotone decreasing on $(N^*-\omega,N^*+\omega)$. Let us define $M:=b(N^*-\omega)-d(N^*-\omega)>0$. Choose $\delta>0$ such that 
$$ \delta \leq \min \left\{\frac{\varepsilon}{3},\frac{\varepsilon}{3}\frac{d}{\gamma}, \frac{M}{d_I},\frac{\omega}{3}\right\}.$$
Assume that the initial data are given such that $|S(0)-N^*|<\delta$, $I(0)<\delta$ and $R(0)=\int_{z_{\min}}^{z_{\max}} r(0,z)dz< \delta$. Then we show that $|S(t)-N^*|<\varepsilon$, $I(t)<\varepsilon$ and $R(t)=||r(t,\cdot)||< \varepsilon$ holds for all $t\geq 0$.

Since $I(t)\leq I(0) e^{-qt} <\delta e^{-qt}$, we have $I(t)<\frac{\varepsilon}{3}<\varepsilon$.
Further, from $R'(t)< \gamma \delta - d R(t),$ we easily find $R(t)<\max \{\delta,\frac{\gamma \delta}{d}\}$ and $R(t)< \frac{\varepsilon}{3}<\varepsilon$ is also guaranteed. Now we consider the susceptible population. Note that $S(t)\leq N(t)$. Let $\hat \nu(t)$ be the solution of the comparison equation to \eqref{eq:ODE_N},
$\hat \nu'(t)=b(\hat \nu(t))-d\hat \nu(t)$. Then for given initial value $\hat \nu(0) \in (N^*-3\delta,N^*+3\delta)$, the solution $\hat \nu(t)$  converges to $N^*$ monotonically. Hence,
we have that ${N(0)=S(0)+I(0)+R(0)<N^*+3\delta}$ implies $S(t)<N^*+3\delta<N^*+\varepsilon$ for all $t>0$.
It remains to show that $S(t)\geq N^*-\varepsilon$ for $t\geq 0$. From the assumptions on $b(N)$, there is a value $\bar N \in (N^*-\omega,N^*)$ such that $b(\bar N) - d \bar N- \delta d_I=0$. Recall that $$N'(t)>b(N(t))-d N(t)- d_I \delta.$$
Let $\breve \nu (t)$ be the solution of the comparison equation $\breve \nu'(t)=b(\breve\nu(t))-d\breve\nu(t)- d_I \delta$, with $\breve\nu(0)=N^*-\omega$. Since
$N(0)> N^*-3\delta>N^*-\omega$, and $\breve\nu(t)$ converges monotonically to $\bar N$, from $N(t)>\breve\nu(t)$ we find that $N(t)> N^*-\omega$. 
Observe that $N(t)=S(t)+I(t)+R(t)$ guarantees that $S(t)>N(t)-\frac{2\varepsilon}{3}$, since $I(t)<\frac{\varepsilon}{3}$ and $R(t)<\frac{\varepsilon}{3}$ hold. In particular,  $S(t)>N^*-\omega-\frac{2\varepsilon}{3}> N^*-\varepsilon$.

Therefore we can conclude the global asymptotic stability of the disease free equilibrium in $X$.
 \qed\\
\ \\
Define the threshold quantity $\tilde R_0:=\frac{\beta}{\gamma+d}$, that coincides with the basic reproduction number in case of non-lethal diseases.
\begin{theorem}[Global stability of DFE]
If $\tilde R_0<1$, the DFE is globally asymptotically stable.  
\end{theorem}
\textit{Proof.} First we show that the components $I$ and $R$ of \eqref{sys:mod1_PDE} converge to zero for $t \to \infty$. Let $q:=\beta \left(\frac{1}{R_0}-1\right)>0$. Then from \eqref{sys:mod1_PDE}, for any solution $I(t)>0$ we have 
$$I'(t)=\beta I(t)\left(\frac{S(t)}{N(t)}-\frac{1}{R_0}\right) < - q I(t),$$
hence $I(t)<I(0) e^{-qt}$ for all $t>0$. Further, from \eqref{dotRt_rel_1}, we have the estimate
$$R'(t)< \gamma I(0) e^{-qt} - d R(t),$$ which implies $R(t)\to 0$ as $t \to \infty$.\\
\ \\
Next we prove that there is an $\eta>0$ such that for any positive solution, $N_\infty:=\liminf_{t\to \infty} N(t) >\eta$.
Observe that the total population $N=S+I+R$ satisfies the equation
\begin{equation}
N'(t)=b(N(t))-dN(t)-d_I I(t).
\label{eq:ODE_N}
\end{equation}
\noindent Recall the Assumption \ref{ass:bN} on $b(N)$ and $d$. Fix a value $\xi<1$ such that $b'(0)\xi>d$. Then we can choose a small $\eta>0$ such that $b(N)>b'(0)\xi N$ for $N \in (0,\eta)$.\\
\indent Now we show that $N_\infty \geq \eta$. Assume the contrary. Then, there are two possibilities:\\
(i) there is a $t_*$ such that $N(t)<\eta$ for all $t>t_*$, or\\
(ii) there is a sequence $t_k \to \infty$ as $k \to \infty$ such that $N(t_k)=\eta$ and $N'(t_k)\leq 0$.\\
\ \\
If (i) holds, then for $t>t_*$ we have from \eqref{eq:ODE_N}
$$N'(t) > (b'(0)\xi-d)N(t)-d_I I(0) e^{-qt}.$$
Consider the ODE
\begin{equation}
x'(t)=Ax(t)-Be^{-Ct},
\label{eq:x_compsys}
\end{equation}
with $A,\,B,\,C>0$, which has the solution
$$x(t)=\frac{B e^{-C t}}{A+C}+e^{A t} \left( x(0)- \frac{B}{A+C}\right).$$
Then clearly $x(t)\to \infty$ whenever $\left( x(0)- \frac{B}{A+C}\right)>0$. 
We can use the equation \eqref{eq:x_compsys} as a comparison system to \eqref{eq:ODE_N} for $t>t_*$ with $x(0)=N(0)$, and $ A=(b'(0)\xi-d)>0, B=d_I I(0)>0, C=q>0$. From $\tilde R_0<1$ it follows that
$d_I<d_I+d+\gamma-\beta=q$. Then $\frac{B}{A+C}<\frac{B}{C}=\frac{d_I}{q}I(0)<I(0)\leq N(0)$, hence $N(t)\to \infty$ and we have
found a contradiction to $N(t) < \eta$.\\
\ \\
On the other side, if (ii) holds, then from $I(t_k) \to 0$ and $b(\eta)-d\eta>0$, there is a $k_*$ such that $d_I I(t_k)< b(\eta)-d\eta$ for $k>k_*$. But then
$$0\geq N'(t_k)=b(N(t_k))-dN(t_k)-d_I I(t_k) > b(\eta)-d\eta-d_I I(t_k)>0,$$
which is a contradiction. We conclude that $N_\infty \geq \eta$. From the fluctuation lemma, there is a sequence $t_k \to \infty$ such that $N(t_k) \to N_\infty$ and $N'(t_k)\to 0$.
Now, considering equation \eqref{eq:ODE_N} with $t=t_k$, and taking the limit for $k\to \infty$, we obtain 
$$0=b(N_\infty)-dN_\infty,$$
which has the only solutions $N_\infty=0$ and $N_\infty=N^*$. Clearly, only $N_\infty=N^*$ is possible. By the same argument, the fluctuation lemma provides us with $N^\infty=N^*$, thus $N(t)\to N^*$. Given that $I(t)\to 0$ and $R(t)\to 0$, we proved that $S(t)\to N^*$. Thus every positive solution converges to the DFE. Moreover, from Theorem \ref{thmDFElocstab}, $R_0<\tilde R_0<1$ guarantees the local stability of the DFE, and the proof is complete.
 \qed\pagebreak

\begin{prop}[Instability of DFE]
If $R_0>1$ the DFE is unstable.
\end{prop}
\textit{Proof.} By straightforward computation, the linearization of the $I$-equation of model (M1) about the DFE yields the simple equation
\begin{equation*}
\label{sys:linsysP0_M1_2}
v'(t) =(\beta-\gamma-d-d_I)v(t)=(\gamma+d+d_I)(R_0-1)v(t),
\end{equation*}
where $v(t)$ is the linear perturbation of $I_*=0$. It is now evident that if $R_0>1$, the DFE repels the $I$ component thus it is unstable. \qed

\section{Connection with ODE models}
\label{sec:connODEs}
In Sect.~\ref{sec:intro} we mentioned several previous works which suggest compartment models for waning immunity and immune system boosting in terms of ODE systems \citep{Grenfell2012,Lavine2011,Heffernan2009,Glass2003b,Dafilis2012}. Thanks to the well-known method of lines we can obtain such ODE systems from our model (M1).\\
\indent The method of lines, mostly used to solve parabolic PDEs, is a technique in which all but one dimension are discretized \cite[see, e.~g.][]{MOLbook}. In our case, we shall discretize the level of immunity ($z$) and obtain a system of ordinary differential equations in the time variable.\\
\ \\
Let us define a sequence $\left\{z_j\right\}_{j\in \N}$, with $h_j:=z_{j+1}-z_j>0$, for all $j \in \N$. To keep the demonstration as simple as possible, we choose a grid with only few points, $z_1:=\zm<\zw<\zf<\zM$ and assume that $h_j=1$ for all $j$. We define the following three subclasses of the R population:
\begin{itemize}\itemsep0.3cm
	\item $R_F(t):= r(t,\zf)$ immune hosts with high level of immunity at time $t$. As their immunity level is quite high, these individuals do not experience immune system boosting. Level of immunity decays at rate $\mu:=g(\zf)>0$.
	\item $R_W(t):= r(t,\zw)$ immune hosts with intermediate level of immunity at time $t$. These individuals can get immune system boosting and move to $R_F$.  Level of immunity decays at rate $\nu:=g(\zw)>0$.
	\item $R_C(t):= r(t,\zm)$ immune hosts with critically low level of immunity at time $t$. With probability $\theta$ boosting moves $R_C$ individuals to $R_W$ (respectively, with probability $(1-\theta)$ to $R_F$).  Level of immunity decays at rate $\sigma:=g(\zm)>0$. If they do not get immune system boosting, these hosts become susceptible again. 
	\end{itemize}

\noindent In the following we show how to proceed in absence of immune system boosting. When immune system boosting is added, the approximation technique remains unchanged.\\
\ \\
For practical reasons we write the equation for $r(t,z)$ in model (M1) in the form
\begin{equation}
\label{eq:PDEr_MOL_noboost}
\frac{\partial}{\partial t} r(t,z) = \frac{\partial}{\partial z} \bigl(g(z) r(t,z)\bigr) -d r(t,z), 
\end{equation}
with boundary condition $R_{\zM}(t):=r(t,\zM)=\gamma I(t) / g(\zM)$.

\noindent Using forward approximation for the $z$-derivative in \eqref{eq:PDEr_MOL_noboost}, we obtain, e.~g. for $R_F(t)$,
\begin{align*}
{R_F}'(t) & =\frac{\partial}{\partial t} r(t,\zf)\\
          &  = \frac{\partial}{\partial z} \bigl(g(\zf) r(t,\zf)\bigr) -d r(t,\zf)\\
					& \approx \frac{g(\zM) r(t,\zM)- g(\zf) r(t,\zf)}{\underbrace{\zM-\zf}_{=1}} -d r(t,\zf)\\
					& = g(\zM) R_{\zM}(t)- \mu R_F(t) -d R_F(t)\\
					& = \gamma I(t)- (\mu+d) R_F(t).
\end{align*}
Analogously we find equations for $R_W$ and $R_C$. Altogether we obtain a system in which a linear chain of ODEs replaces the PDE for the immune class,
\begin{align*}
S'(t) & = b(N(t)) -\beta \frac{S(t)I(t)}{N(t)}-dS(t)+\sigma R_C(t),\\[0.2em]
I'(t) & = \beta \frac{S(t)I(t)}{N(t)}-(\gamma+d+d_I)I(t),\\[0.2em]
{R_F}'(t) & = \gamma I(t)-\mu R_F(t) -dR_F(t),\\[0.2em]
{R_W}'(t) & =\mu R_F(t)-\nu R_W(t)-d R_W(t),\\[0.2em]
{R_C}'(t) & = \nu R_W(t)-\sigma R_C(t) -d R_C(t).
\end{align*}
The linear chain of ODEs provides a rough approximation of the PDE in model (M1). Indeed, with the method of lines we approximate the PDE dynamics considering changes only at the grid points ($\zm,\;\zw,\;\zf,\;\zM$), whereas the dynamics remains unchanged in the interval between one grid point and the next one. In other words, we average out over the immunity level in one immunity interval $[z_{j},z_{j+1}]$ and consider as representative point of the interval the lowest boundary $z_j$. This is also the reason why we do not have a differential equation for $R_{\zM}(t)$, but simply a boundary condition at this point.\\
\ \\
\noindent \textbf{Including immune system boosting}\\
Now we use the same approximation scheme for the full model (M1), which includes also immune system boosting to any higher immunity level. To this purpose it is necessary to specify the boosting probability $p(z,\tilde z)$ with $z \in [\zm,\zM]$, which we choose as follows
\begin{displaymath}
p(z,\tilde z) =
\begin{cases}
1 & \text{if } \tilde z=\zw,\text{ and } z=\zf, \\
\theta & \text{if } \tilde z=\zm,\text{ and } z=\zw, \\
1-\theta & \text{if } \tilde z=\zm,\text{ and } z=\zf, \\
0 & \text{otherwise}.
\end{cases}
\end{displaymath}
The integral term in model (M1) is discretized by the mean of a finite sum and the resulting ODE system is 
\begin{equation}
\begin{aligned}
S'(t) & = b(N(t)) -\beta \frac{S(t)I(t)}{N(t)}-dS(t)+\sigma R_C(t),\\[0.3em]
I'(t) & = \beta \frac{S(t)I(t)}{N(t)}-(\gamma+d+d_I)I(t),\\[0.3em]
{R_F}'(t) & = \gamma I(t)-\mu R_F(t) -dR_F(t)+\beta \frac{I(t)}{N(t)}\biggl((1-\theta)R_C(t) +R_W(t)\biggr),\\[0.3em]
{R_W}'(t) & =\mu R_F(t)-\nu R_W(t)-d R_W(t)+\beta \frac{I(t)}{N(t)}\biggl(\theta R_C(t)- R_W(t)\biggr),\\[0.3em]
{R_C}'(t) & = \nu R_W(t)-\sigma R_C(t) -d R_C(t)- \beta\frac{I(t)}{N(t)}R_C(t).
\end{aligned}
\label{sys:mod1}
\end{equation}
To obtain an ODE system from model (M2), that is, the special case in which boosting always restores the maximal level of immunity, we include all hosts who get immune system boosting in the boundary condition,
\begin{equation*}
R_{\zM}(t)=\frac{1}{g(\zM)}\left(\gamma I(t) + \beta \frac{I(t)}{N(t)}\left(R_W(t)+R_C(t)\right)\right).
\end{equation*}
This correspond to setting $\theta =0$ in \eqref{sys:mod1}.

\begin{figure}[!]
\centering
\includegraphics[width=0.7\columnwidth]{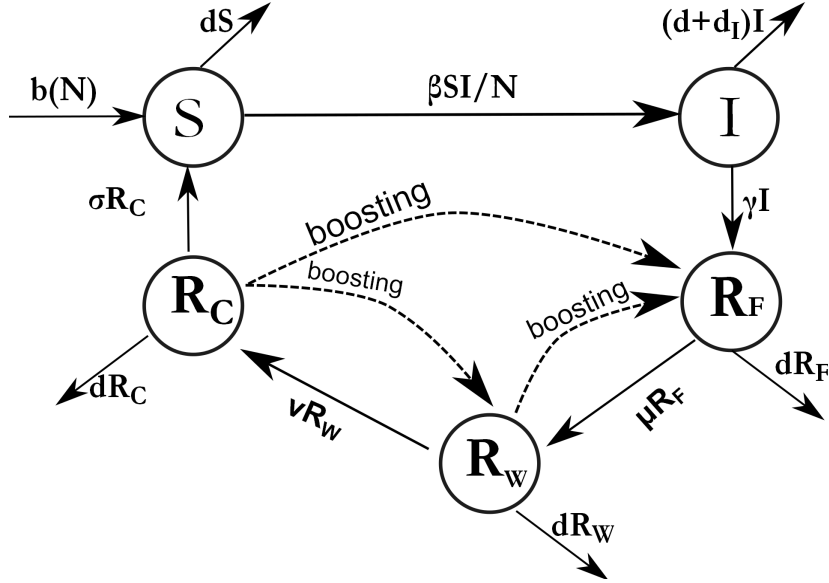}
\caption{Diagram of model \eqref{sys:mod1}. After recovery, individuals enter the $R_F$ class and are protected for a while. Due to contact with infectives, the immune system of partially immune individuals ($R_W$ and $R_C$) can be boosted to a higher level of immunity. Natural decay of the immune status moves individuals from $R_F$ to $R_W$, from $R_W$ to $R_C$ and finally from $R_C$ to $S$.}
\label{Fig:schema_model1}
\end{figure}

\section{Connection with DDE models}
\label{sec:connDDEs}
Delay models with constant delay can be obtained from special cases of model (M1) defining the delay $\tau$ as the duration of immunity induced by the natural infection. We show here how the classical SIRS model with delay, as well as a new SIS model with delay arise from model (M1). 

\subsection{SIRS with constant delay}
Consider the general model (M1) and neglect the boosting effects after recovery, that is, set $c_{max}(z)\equiv 0 \equiv c_0(z)$ for all $z \in [\zm,\zM]$. This means that equations \eqref{eq:mod1_pde_R}--\eqref{eq:BC_mod1_pde_R} change into
\begin{align*}
\frac{\partial }{\partial t}r(t,z) - \frac{\partial }{\partial z}\left(g(z)r(t,z)\right) & = -dr(t,z)\\
g(\zM)r(t,\zM)& = \gamma I(t).
\end{align*}
From our assumptions, the disease-induced immunity lasts for a fix time, $\tau>0$ years, given by
$$ \int_{\zm}^{\zM}\frac{1}{g(x)}\,dx=\tau.$$
We can express the total immune population at time $t$ as the number of individuals who recovered in the time interval $[t-\tau,t]$,
\begin{equation*}
R(t) = \gamma \int_{t-\tau}^{t} I(y)e^{-d(t-y)}\,dy = \gamma \int_0^{\tau} I(t-x)e^{-dx}\,dx.
\label{eq:defR_integr}
\end{equation*}
Differentiation with respect to $t$ yields
\begin{equation}
R'(t) = \gamma I(t)-\gamma I(t-\tau)e^{-d\tau}-dR(t).
\label{eq:ddeR1}
\end{equation}
On the other side, we have the relation \eqref{dotRt_rel_1}.
Comparison between \eqref{eq:ddeR1} and \eqref{dotRt_rel_1} yields
\begin{equation*}
g(\zm)r(t,\zm)=\gamma I(t-\tau)e^{-d\tau},
\label{eq:ddeR_relations}
\end{equation*}
which confirms the model assumptions, namely, $\tau$ time after recovery immune hosts who did not die in the $\tau$ interval of time become susceptible again. 
In other words, if an individual who recovers at time $t_1$ survives up to time $t_1+\tau$, he exits the $R$ class and enters $S$. 
In turn, we find a delay term in the equation for $S$ too, and have a classical SIRS model with constant delay 
\begin{align*}
S'(t) & = b(N(t)) -\beta \frac{S(t)I(t)}{N(t)}-dS(t)+\gamma I(t-\tau)e^{-d\tau},\\[0.3em]
I'(t) & = \beta \frac{S(t)I(t)}{N(t)}-(\gamma+d+d_I)I(t),\\[0.3em]
R'(t) & =\gamma I(t)-\gamma I(t-\tau)e^{-d\tau} -dR(t).
\end{align*}
This model was studied by \citet{Taylor2009}.

\subsection{A new class of SIRS models with constant delay}
\label{sec:newsis}
Consider the special case (M2) in which boosting restores the maximal level of immunity. 
As in the previous case, let $\tau>0$ be the duration of immunity induced by the natural infection. 
With the definition of the characteristic curve $\zeta(t)=\varphi(t;0,\zM)$ (cf. Sect.~\ref{sec:ex_maxlevel_boostingmodel}), we have the relation
\begin{displaymath}
\zm=\varphi(t;t-\tau,\zM).
\end{displaymath}
Thus, solving along the characteristics, for $t>\tau$ we have
\begin{align*}
r(t,\zm) & = \frac{B(t-\tau)}{g(\zm)} \exp \left( -\int_{t-\tau}^{t}\mu(u,\varphi(u;t-\tau,\zM))\,du\right)\\
         & = \frac{B(t-\tau)}{g(\zm)} 
		\exp \left( -\int_{t-\tau}^{t}\left(d+\beta\frac{I(u)}{N(u)}\right)\,du\right)\\
         & = \frac{I(t-\tau)}{g(\zm)} \left(\gamma +  \frac{\beta R(t-\tau)}{N(t-\tau)} \right)
     \exp \left( -d\tau -\beta\int_{t-\tau}^{t}\frac{I(u)}{N(u)}\,du\right).
\end{align*}
Assume for simplicity that there is no disease-induced death ($d_I=0$) and that the total population is constant $N(t)\equiv 1$. 
Then, $R(t)=1-S(t)-I(t)$ and we end up with the system
\begin{align*}
S'(t) & = d(1-S(t)) -\beta I(t)S(t)\\
& \quad +I(t-\tau) \left(\gamma + \beta \left(1-S(t-\tau)-I(t-\tau)\right)\right)
     \exp \left( -d\tau - \beta\int_{t-\tau}^{t}I(u)\,du\right),\\[0.3em]
I'(t) & = \beta I(t)S(t) -(\gamma+d) I(t).
\end{align*}
Immunity loss occurs at time $t$ either for hosts who recovered from infection at time $t-\tau$, 
or for host who, being immune, were exposed to the pathogen at time $t-\tau$, and in the interval of time $[t-\tau,t]$ did neither die nor come in contact with infectives. It seems that such a SIS model with constant delay has not been presented in previous literature.

\section{Discussion}
\label{sec:discussion}
This paper provides a general framework for modeling waning immunity and immune system boosting. The model (M1) combines the in-host perspective with the population dynamics, while keeping the number of equations as low as possible. To the best of our knowledge, only the model proposed by \citet{Martcheva2006} achieved a similar result. However \citet{Martcheva2006} give a different interpretation of immune system boosting (observed only during the infectious period).\\
\ \\
Although, when first examined, the immune hosts equation \eqref{eq:mod1_pde_R} in model (M1) could look like a size-structured model, we explained in Sect.~\ref{sec:model} that this is not the case. Equation \eqref{eq:mod1_pde_R} describes a physiologically structured population whose dynamics includes a transport process and jumps. Such kind of models are very rare in mathematical biology, possibly because they are not easily derived or because their qualitative analysis turns out to be very challenging. A somehow similar example for the transmission of a microparasite was suggested in \cite{White1998}. Their model is given by two PDEs for uninfected and infected hosts structured by the level of immunity. However, existence/uniqueness, positivity of solutions, or stability analysis were not considered in \cite{White1998}. In contrast, we provided in Sect.~\ref{sec:exiuni} results on existence of a unique classical solution and on positive invariance, whereas in Sect.~\ref{sec:dfe} we discussed the global behavior of solutions with respect to the disease-free equilibrium.\\
\ \\
\noindent We have shown that (M1) is a general framework which allows to recover previous models from the literature, such as systems of ODEs (Sect.~\ref{sec:connODEs}), as well as models with constant delay (Sect.~\ref{sec:connDDEs}). Systems of ODEs presented, e.~g. in \cite{Grenfell2012,Lavine2011,Heffernan2009} can be obtained by discretizing the PDE for the immune population in model (M1) with the help of the method of lines. This method, however, provides only a rough approximation of the structured population.\\
\indent The system suggested in \cite{Dafilis2012} includes a boosting term where a constant value $\nu>0$, the boosting rate, is multiplied by the force of infection $\beta$. Our work includes the case $\nu\leq 1$, nevertheless our results can be easily extended to the case of any positive $\nu$.\\
\indent Models with constant delay, such as those in \cite{Taylor2009,Kyrychko2005}, are obtained from model (M1) neglecting immune system boosting. The delay represents the duration of the immunity after natural infection. Moreover we could show that there is a new and so far not studied class of SIS models with constant delay which can be obtained from model (M2), the special case in which boosting always restores the maximal level of immunity. Both limit cases, no boosting or boosting always to the maximal level of immunity, can be reduced to systems of equations with constant delays. The same cannot be stated for the more general case in which boosting can restore any higher immunity level. In this case, indeed, the obtained delay is of state-dependent type. As the connection between model (M1) and system of equations with state-dependent delay is nontrivial, we leave it as a future project.\\
\ \\
Setting up our model we did not restrict ourselves to a particular pathogen. Our goal is to provide a general framework for immune response, immune system boosting and waning immunity in the context of population dynamics. Choosing appropriate model coefficients, which can be deduced from experimental data, the system (M1) can be adapted to model epidemic outbreaks, such as measles, chickenpox, rubella, pertussis. In this context it is worth mentioning that, to date, not many data are available about immune system boosting after natural infection. Recent experimental data provide some information on the effects of vaccines on the antibody level, see for example \cite{Luo2012} for hepatitis B in mice or \cite{Amanna2007,Li2013} for measles, mumps, rubella and influenza in humans.\\
\ \\
We conclude with some remarks on vaccination. As we have mentioned in Sect.~\ref{sec:intro}, immunization is not only the result of natural infection, but also of vaccination and transmission of maternal antibodies (passive-immunity). In a highly vaccinated population there are a lot of individuals with vaccine-induced immunity and few infection cases, as well as more individuals with low level of
immunity. In other words, if a high portion of the population gets the vaccine, there are very few chances for exposure to the pathogen and consequently for immune system boosting in protected individuals. This might be one aspect which causes recurrent outbreaks of, e.~g. measles or pertussis in highly vaccinated populations.\\
\indent In order to understand how in-host processes like waning immunity and immune system boosting are related to the dynamics of the population, and to provide a correct and whole-comprehensive mathematical formulation for these phenomena we chose not to include vaccination in our model (M1). 
One natural extension of our work in the future is to include vaccine-induced immunity in our framework. The resulting approach may provide a general setting which connects various models in the literature, for example in \cite{Mossong1999,Glass2003b,Grenfell2012,Arino2006,Mossong2003}.

\begin{acknowledgements}
MVB was supported by the ERC Starting Grant No. 259559 as well as by the European Union and
the State of Hungary, co-financed by the European Social Fund in the framework of T\'AMOP-4.2.4.
A/2-11-1-2012-0001 National Excellence Program. GR was supported by Hungarian Scientific Research Fund OTKA K109782 and T\'AMOP-4.2.2.A-11/1/KONV-2012-0073 "Telemedicine focused
research activities on the field of Mathematics, Informatics and Medical
sciences".
\end{acknowledgements}

\newpage
\appendix
\section{Proof of Theorem \ref{thm:existence_M1_solux}}
In the following we show the continuous differentiability of the map $Q$, which is necessary to have existence and uniqueness of a classical solution of the abstract Cauchy problem \eqref{abst_cauch_probl_X}. Continuous differentiability of $Q$ can be shown in two steps:
(a) First we determine the existence of the operator $DQ(x;w)$, for all $x,\,w \in X$, defined by  
\begin{displaymath}
DQ(x;w):=\lim \limits_{h \to 0}\frac{Q(x+hw)-Q(x)}{h}.
\end{displaymath}
(b) Second we show that the operator $DQ(x;\cdot)$ is continuous in $x$, that is
\begin{displaymath}
\lim \|DQ(x;\cdot) -DQ(y;\cdot)\|_{OP}=0 \qquad \mbox{for} \qquad \normx{x-y}\to 0,
\end{displaymath}
where $\|\cdot\|_{OP}$ is the operator norm.\\
\ \\
For simplicity of notation we write
\begin{align*}
DQ_1(x; w) &: = P_1(x;w)-P_2(x;w),\\[0.1em]
DQ_2(x; w) &: = P_2(x;w),\\[0.1em]
DQ_3(x; w) &: = -P_3(x;w)+P_4(x;w), 
\end{align*}
where
\begin{align*}
P_1(x;w) & :=\,\lim \limits_{h \to 0}\frac{b(\hat x+h \hat w)-b(\hat x)}{h},\\[0.2em]
P_2(x;w) & :=\,\lim \limits_{h \to 0}\frac{1}{h} \beta \left( \frac{(x_1+hw_1)(x_2+hw_2)}{\hat x+h \hat w}-\frac{x_1x_2}{\hat x}\right),\\[0.2em]
P_3(x;w)& :=\,\lim \limits_{h \to 0}\frac{1}{h} \beta \left( 
\frac{(x_2+hw_2)(x_3(z)+hw_3(z))}{\hat x+h\hat w} - \frac{x_2x_3(z)}{\hat x}\right),\\[0.2em]
 P_4(x;w) & :=\,\lim \limits_{h \to 0}\frac{1}{h} \beta \left( 
\frac{(x_2+hw_2)\int_{\zm}^{z} (x_3(u)+hw_3(u))\,p(z,u)\,du}{\hat x+h\hat w} - \frac{x_2\int_{\zm}^{z} x_3(u)p(z,u)\,du}{\hat x}\right).
\end{align*}

\noindent \textit{Proof of (a).} We compute the limit for the first component of $Q_1$,
\begin{align*}
P_1(x;w)& = \lim \limits_{h \to 0}\frac{1}{h} \biggl( b\left(\hat x+h \hat w\right)
-b\left(\hat x+hw_1+hw_2\right)\biggr)\\
& \phantom{=} + \lim \limits_{h \to 0}\frac{1}{h} \biggl( b\left(\hat x+hw_1+hw_2\right)
-b\left(\hat x+hw_1\right)\biggr)\\
& \phantom{=}+\lim \limits_{h \to 0}\frac{1}{h} \bigl( b\left(\hat x+hw_1\right)
-b\left(\hat x\right)\bigr)\\[0.3em]
& = \lim \limits_{h \to 0}\frac{b\left(\hat x+h \hat w\right)
-b\left(\hat x+hw_1+hw_2\right)}{h \int_{\zm}^{\zM}w_3(u)\,du} \int_{\zm}^{\zM}w_3(u)\,du\\
& \phantom{=} +\lim \limits_{h \to 0}\frac{b\left(\hat x+hw_1+hw_2\right)
-b\left(\hat x+hw_1\right)}{h w_2}\,w_2\\
& \phantom{=} + \lim \limits_{h \to 0}\frac{b\left(\hat x+hw_1\right)
-b\left(\hat x\right)}{h\,w_1} \,w_1\\[0.3em]
& = b'(\hat x)  \hat w.
\end{align*}

\noindent For the second term in $Q_1(x)$ we have
\begin{align*}
P_2(x;w)
& = \beta \lim \limits_{h \to 0}\frac{1}{h}  \left( \frac{(x_1+hw_1)(x_2+hw_2)}{\hat x+h \hat w}-
\frac{(x_1+hw_1)(x_2+hw_2)}{\hat x + h w_2+h\int_{\zm}^{\zM}w_3(u)\,du}\right)\\
& \phantom{=} + \beta \lim \limits_{h \to 0}\frac{1}{h} \left(\frac{(x_1+hw_1)(x_2+hw_2)}{\hat x + h w_2+h\int_{\zm}^{\zM}w_3(u)\,du}
- \frac{x_1(x_2+hw_2)}{\hat x + h w_2+h\int_{\zm}^{\zM}w_3(u)\,du}\right)\\
& \phantom{=} + \beta \lim \limits_{h \to 0}\frac{1}{h} \left(\frac{x_1(x_2+hw_2)}{\hat x + h w_2+h\int_{\zm}^{\zM}w_3(u)\,du}
- \frac{x_1(x_2+hw_2)}{\hat x +h\int_{\zm}^{\zM}w_3(u)\,du}\right)\\
& \phantom{=} + \beta \lim \limits_{h \to 0}\frac{1}{h} \left(\frac{x_1(x_2+hw_2)}{\hat x +h\int_{\zm}^{\zM}w_3(u)\,du}
-\frac{x_1x_2}{\hat x +h\int_{\zm}^{\zM}w_3(u)\,du}\right)\\
& \phantom{=} +  \beta \lim \limits_{h \to 0}\frac{1}{h} \left( \frac{x_1x_2}{\hat x +h\int_{\zm}^{\zM}w_3(u)\,du}
-\frac{x_1x_2}{\hat x}\right)\\[0.3em]
& = -\beta \lim \limits_{h \to 0} \frac{(x_1+hw_1)(x_2+hw_2)}{\left(\hat x+h \hat w\right)\left(\hat x + h w_2+h\int_{\zm}^{\zM}w_3(u)\,du\right)}\,w_1\\
& \phantom{=} + \beta \lim \limits_{h \to 0}  \frac{x_2+hw_2}{\hat x + h w_2+h\int_{\zm}^{\zM}w_3(u)\,du}\,w_1\\
& \phantom{=} - \beta \lim \limits_{h \to 0}  \frac{x_1(x_2+hw_2)}{\left(\hat x + h w_2+h\int_{\zm}^{\zM}w_3(u)\,du\right)\left(\hat x +h\int_{\zm}^{\zM}w_3(u)\,du\right)}\,w_2\\
& \phantom{=} + \beta \lim \limits_{h \to 0}  \frac{x_1}{\hat x +h\int_{\zm}^{\zM}w_3(u)\,du}\,w_2\\
& \phantom{=} -  \beta \lim \limits_{h \to 0} \frac{x_1x_2}{\left(\hat x +h\int_{\zm}^{\zM}w_3(u)\,du\right) \hat x}\,\int_{\zm}^{\zM}w_3(u)\,du\\[0.5em]
& = \beta \left[ \frac{x_2(\hat x-x_1)}{{\hat x}^2}\,w_1
 + \frac{x_1(\hat x -x_2)}{{\hat x}^2} \,w_2
 -\frac{x_1x_2}{{\hat x}^2}\,\int_{\zm}^{\zM}w_3(u)\,du\right].
\end{align*}
The first term in $Q_3(x)$:
\begin{align*}
P_3(x;w) & = \beta \lim \limits_{h \to 0}\frac{1}{h} \left( \frac{(x_2+hw_2)(x_3(z)+hw_3(z))}{\hat x+h\hat w}-
\frac{(x_2+hw_2)(x_3(z)+hw_3(z))}{\hat x+h w_2+h\int_{\zm}^{\zM}w_3(u)\,du}\right)\\
&\phantom{=} + \beta \lim \limits_{h \to 0}\frac{1}{h} \left( \frac{(x_2+hw_2)(x_3(z)+hw_3(z))}{\hat x+h w_2+h\int_{\zm}^{\zM}w_3(u)\,du}
 - \frac{x_2(x_3(z)+hw_3(z))}{\hat x+h w_2+h\int_{\zm}^{\zM}w_3(u)\,du}\right)\\
& \phantom{=}+ \beta \lim \limits_{h \to 0}\frac{1}{h} \left(  \frac{x_2(x_3(z)+hw_3(z))}{\hat x+ w_2+h\int_{\zm}^{\zM}w_3(u)\,du}
- \frac{x_2(x_3(z)+hw_3(z))}{\hat x+h\int_{\zm}^{\zM}w_3(u)\,du}\right)\\
&\phantom{=} + \beta \lim \limits_{h \to 0}\frac{1}{h} \left( \frac{x_2(x_3(z)+hw_3(z))}{\hat x+h\int_{\zm}^{\zM}w_3(u)\,du}
- \frac{x_2x_3(z)}{\hat x+h\int_{\zm}^{\zM}w_3(u)\,du}\right)\\
& \phantom{=}+ \beta \lim \limits_{h \to 0}\frac{1}{h} \left( \frac{x_2x_3(z)}{\hat x+h\int_{\zm}^{\zM}w_3(u)\,du}
- \frac{x_2x_3(z)}{\hat x}\right).
\end{align*}

\noindent Hence we have
\begin{align*}
P_3(x;w)  & = -\beta \lim \limits_{h \to 0} \frac{(x_2+hw_2)(x_3(z)+hw_3(z))}{\left(\hat x+h\hat w\right)\left(\hat x+h w_2+h\int_{\zm}^{\zM}w_3(u)\,du\right)}\,w_1\\
& \phantom{=}+ \beta \lim \limits_{h \to 0}\frac{x_3(z)+hw_3(z)}{\hat x+h w_2+h\int_{\zm}^{\zM}w_3(u)\,du}\,w_2\\
& \phantom{=} -\beta \lim \limits_{h \to 0}  \frac{x_2(x_3(z)+hw_3(z))}{\left(\hat x+h w_2+h\int_{\zm}^{\zM}w_3(u)\,du\right)\left(
\hat x+h\int_{\zm}^{\zM}w_3(u)\,du\right)}\,w_2\\
& \phantom{=} + \beta \lim \limits_{h \to 0}\frac{x_2}{\hat x+h\int_{\zm}^{\zM}w_3(u)\,du}\,w_3(z)\\
& \phantom{=}   - \beta \lim \limits_{h \to 0}\frac{x_2x_3(z)}{\hat x\left(\hat x+h\int_{\zm}^{\zM}w_3(u)\,du\right)}\,\int_{\zm}^{\zM}w_3(u)\,du\\[0.3em]
& = \beta \biggl[-\frac{x_2x_3(z)}{{\hat x}^2}\,w_1+\frac{x_3(z)(\hat x -x_2)}{{\hat x}^2}\,w_2
+ \frac{x_2}{\hat x}\,w_3(z) - \frac{x_2x_3(z)}{{\hat x}^2}\,\int_{\zm}^{\zM}w_3(u)\,du\biggr].
\end{align*}
Analogously, compute the last term in $Q_3(x)$,
\begin{align*}
P_4(x;w) 
& = \beta \lim \limits_{h \to 0}\frac{\int_{\zm}^{z} (x_3(u)+hw_3(u))\,p(z,u)\,du}{\hat x+h\hat w} \,w_2\\
& \phantom{=} - \beta \lim \limits_{h \to 0}
\frac{x_2\int_{\zm}^{z} (x_3(u)+hw_3(u))\,p(z,u)\,du}{\left(\hat x+h\hat w\right)\left(\hat x+h w_2+ h \int_{\zm}^{\zM}w_3(u)\,du\right)}\,w_1\\
& \phantom{=} + \beta \lim \limits_{h \to 0}
\frac{x_2}{\hat x+h w_2+ h \int_{\zm}^{\zM}w_3(u)\,du}\,\int_{\zm}^{z} w_3(u)\,p(z,u)\,du\\
& \phantom{=} -\beta \lim \limits_{h \to 0}\frac{x_2\int_{\zm}^{z} x_3(u)\,p(z,u)\,du}{\left(\hat x+h w_2+ h \int_{\zm}^{\zM}w_3(u)\,du\right)
\left(\hat x+ h \int_{\zm}^{\zM}w_3(u)\,du\right)}\,w_2\\
& \phantom{=} - \beta \lim \limits_{h \to 0} \frac{x_2\int_{\zm}^{z} x_3(u)\,p(z,u)\,du}{\left(\hat x+ h \int_{\zm}^{\zM}w_3(u)\,du\right)\hat x}\,\int_{\zm}^{\zM}w_3(u)\,du \\[0.3em]
& = \beta \left[ - \frac{x_2\int_{\zm}^{z} x_3(u)\,p(z,u)\,du}{{\hat x}^2}\,w_1
+ \frac{(\hat x -x_2)\int_{\zm}^{z} x_3(u)\,p(z,u)\,du}{{\hat x}^2}\,w_2\right.\\
& \qquad \quad \left.+ \frac{x_2\left(\hat x -\int_{\zm}^{z} x_3(u)\,p(z,u)\,du\right)}{{\hat x}^2}\int_{\zm}^{\zM}w_3(u)\,du
 \right].
\end{align*}\pagebreak

\noindent To prove that the operator $DQ(x;\cdot)$ is continuous in $x$ we consider the norm\linebreak
 $\|DQ(x;\cdot) -DQ(y;\cdot) \|_{OP}$, that is
\begin{align*}
\lefteqn{ \underset{ \normx{w}\leq 1 }{\sup} \normx{DQ(x;w) -DQ(y;w)}}\\
& = \underset{ \normx{w}\leq 1 }{\sup} \biggl\{|P_1(x;w)-P_1(y;w)-P_2(x;w)+P_2(y;w)| + |P_2(x;w)-P_2(y;w)|\\
& \qquad \qquad + \int_{\zm}^{\zM} |(P_3(y;w)-P_3(x;w)+P_4(x;w)-P_4(y;w))(z)|\,dz\biggr\},
\end{align*}
and show that 
$$ \|DQ(x;\cdot) -DQ(y;\cdot) \|_{OP} \to 0, \qquad \mbox{for } \quad \normx{x- y}\to 0.$$
We estimate the operator norm as follows
\begin{displaymath}
	\underset{ \normx{w}\leq 1 }{\sup} \normx{DQ(x;w) -DQ(y;w)}\leq \underset{ \normx{w}\leq 1 }{\sup} \sum \limits_{j=1}^{4} T_j(x,y;w),
\end{displaymath}
with
\begin{align*}
T_1(x,y;w) & = |P_1(x;w)-P_1(y;w)|,\\
T_2(x,y;w) & = 2|P_2(x;w)-P_2(y;w)|,\\
T_3(x,y;w) & = \int_{\zm}^{\zM} |(P_3(x;w)-P_3(y;w))(z)|\,dz,\\
T_4(x,y;w) & =\int_{\zm}^{\zM}|(P_4(x;w)-P_4(y;w))(z)|\,dz.
\end{align*}
Then we show the convergence to zero of the above sum. It is obvious that 
$$ \underset{ \normx{w}\leq 1 }{\sup} T_1(x,y;w) \; \to 0, \qquad \mbox{for } \normx{x- y}\to 0, $$
as $b$ is continuously differentiable (see Assumption~\ref{ass:bNdN}),
\begin{align*}
	|P_1(x;w)-P_1(y;w)| & = |b'(\hat x)  \hat w - b'(\hat y)  \hat w|\;\leq |b'(\hat x) - b'(\hat y)| \normx{w}.
\end{align*}
The term in $T_2(x,y;w)$ can be estimated as follows:
\begin{align*}
|P_2(x;w)-P_2(y;w)| & \leq \underbrace{\beta \left| \left( \frac{x_2(\hat x-x_1)}{{\hat x}^2}-\frac{y_2(\hat y-y_1)}{{\hat y}^2}\right)\right|\,|w_1|}_{=:L_1(x,y;w)}\\
& \phantom{\leq} +\underbrace{\beta \left| \left( \frac{x_1(\hat x -x_2)}{{\hat x}^2}-\frac{y_1(\hat y -y_2)}{{\hat y}^2}\right)\right| \,|w_2|}_{=:L_2(x,y;w)}\\
& \phantom{\leq} +\underbrace{\beta \left| \left(\frac{x_1x_2}{{\hat x}^2}-\frac{y_1y_2}{{\hat y}^2}\right)\right|\,\int_{\zm}^{\zM}|w_3(u)|\,du}_{=:L_3(x,y;w)}.
\end{align*}
Since the addends of the last sum are all similar, we show convergence to zero for only one of them.
\begin{align*}
\lefteqn{ \underset{ \normx{w}\leq 1 }{\sup}\;L_3(x,y;w)}\\[0.5em]
& \leq \beta \left| \frac{x_1x_2}{{\hat x}^2}-\frac{y_1y_2}{{\hat y}^2}\right|\\
& \leq \beta \biggl( \left| \frac{x_1x_2}{{\hat x}^2}- \frac{x_1x_2}{\hat x\hat y}\right|+
\left| \frac{x_1x_2}{\hat x \hat y}-\frac{x_1y_2}{\hat x \hat y}\right| +  \left| \frac{x_1y_2}{\hat x \hat y}-\frac{x_1y_2}{{\hat y}^2}\right|
+ \left|\frac{x_1y_2}{{\hat y}^2}-\frac{y_1y_2}{{\hat y}^2}\right|\biggr)\\[0.5em]
& \leq \beta \biggl( \left| \frac{x_1x_2(\hat y-\hat x)}{{\hat x}^2\hat y}\right|+
\left| \frac{x_1(x_2-y_2)}{\hat x \hat y}\right| +  \left| \frac{x_1y_2(\hat y-\hat x)}{\hat x {\hat y}^2}\right|
+ \left|\frac{(x_1-y_1)y_2}{{\hat y}^2}\right|\biggr)\\[0.5em]
& \leq \frac{3\beta}{|\hat y|} \, \normx{x-y}. 
\end{align*}
It works analogously for the terms $L_1(x,y;w)$ and $L_2(x,y;w)$. Hence we have that 
$$\underset{ \normx{w}\leq 1 }{\sup} T_2(x,y;w)  \; \to 0, \qquad \mbox{for } \normx{x- y}\to 0. $$
\ \\
In a similar way one can estimate $T_3(x,y;w)$ and $T_4(x,y;w)$. We show the computation for $T_4$ which is the most challenging of the two, as it includes double integrals.
\begin{align*}
\lefteqn{T_4(x,y;w)  \leq \underbrace{\beta \int_{\zm}^{\zM} \left|\frac{x_2\int_{\zm}^{z} x_3(u)\,p(z,u)\,du}{{\hat x}^2}
-\frac{y_2\int_{\zm}^{z} y_3(u)\,p(z,u)\,du}{{\hat y}^2}\right|\,|w_1|\,dz}_{=:F_1(x,y;w)}}\\[0.3em]
& + \underbrace{\beta \int_{\zm}^{\zM} \left|\frac{(\hat x -x_2)\int_{\zm}^{z} x_3(u)\,p(z,u)\,du}{{\hat x}^2}
-\frac{(\hat y -y_2)\int_{\zm}^{z} y_3(u)\,p(z,u)\,du}{{\hat y}^2}\right|\,|w_2|\,dz}_{=:F_2(x,y;w)}\\[0.3em]
& + \underbrace{\beta \int_{\zm}^{\zM} \left| \frac{x_2}{{\hat x}}
-\frac{y_2}{{\hat y}}\right|\,\int_{\zm}^{\zM}|w_3(u)|\,du\,dz}_{=:F_3(x,y;w)}\\[0.3em]
& + \underbrace{\beta \int_{\zm}^{\zM} \left| \frac{x_2\int_{\zm}^{z} x_3(u)\,p(z,u)\,du}{{\hat x}^2}
-\frac{y_2\int_{\zm}^{z} y_3(u)\,p(z,u)\,du}{{\hat y}^2}\right|\,\int_{\zm}^{\zM}|w_3(u)|\,du\,dz}_{=:F_4(x,y;w)}.
\end{align*}
Before proceeding to the next estimate, it is useful to observe that 
\begin{align*}
\int_{\zm}^{\zM} \int_{\zm}^{z} x_3(v)p(z,v)\,dv\,dz 
& =\int_{\zm}^{\zM} \int_{v}^{\zM} x_3(v)p(z,v)\,dz\,dv\\
& =\int_{\zm}^{\zM} x_3(v) \underbrace{\int_{v}^{\zM} p(z,v)\,dz}_{=1}\,dv\\
& =\int_{\zm}^{\zM} x_3(v)\,dv.
\end{align*}
Let us now consider the last addend in $T_4(x,y;w)$. We have that
\begin{align*}
\underset{ \normx{w}\leq 1 }{\sup}\; F_4(x,y;w) & \leq \beta  \int_{\zm}^{\zM} \left| \frac{x_2\int_{\zm}^{z} x_3(u)\,p(z,u)\,du}{{\hat x}^2}
-\frac{y_2\int_{\zm}^{z} y_3(u)\,p(z,u)\,du}{{\hat y}^2}\right|\,dz\\[0.5em]
& \leq \beta  \int_{\zm}^{\zM} \left| \frac{x_2\int_{\zm}^{z} x_3(u)\,p(z,u)\,du}{{\hat x}^2}
- \frac{x_2\int_{\zm}^{z} x_3(u)\,p(z,u)\,du}{{\hat x}\hat y}\right|\,dz\\
& \phantom{\leq } + \beta  \int_{\zm}^{\zM} \left| \frac{x_2\int_{\zm}^{z} x_3(u)\,p(z,u)\,du}{{\hat x}\hat y}
-\frac{x_2\int_{\zm}^{z} y_3(u)\,p(z,u)\,du}{{\hat x}\hat y}\right|\,dz\\
& \phantom{\leq } + \beta  \int_{\zm}^{\zM} \left|\frac{x_2\int_{\zm}^{z} y_3(u)\,p(z,u)\,du}{{\hat x}\hat y}
-\frac{x_2\int_{\zm}^{z} y_3(u)\,p(z,u)\,du}{{\hat y}^2}\right|\,dz\\
& \phantom{\leq } + \beta  \int_{\zm}^{\zM} \left|\frac{x_2\int_{\zm}^{z} y_3(u)\,p(z,u)\,du}{{\hat y}^2}
- \frac{y_2\int_{\zm}^{z} y_3(u)\,p(z,u)\,du}{{\hat y}^2}\right|\,dz.
\end{align*}
A similar computation as the one for $\underset{ \normx{w}\leq 1 }{\sup}\;L_3(x,y;w)$ yields
\begin{align*}
\underset{ \normx{w}\leq 1 }{\sup}\; F_4(x,y;w) 
 & \leq \beta  \int_{\zm}^{\zM} \int_{\zm}^{z} |x_3(u)\,p(z,u)|\,du\,dz \, \left|\frac{x_2(\hat y-\hat x)}{{\hat x}^2\hat y}\right|\\
& \phantom{\leq } + \beta  \int_{\zm}^{\zM} \int_{\zm}^{z} |x_3(u)\,p(z,u)|\,du\,dz\, \left|\frac{x_2(\hat y-\hat x)}{\hat x {\hat y}^2}\right|\\
& \phantom{\leq } + \beta  \int_{\zm}^{\zM} \int_{\zm}^{z} |(x_3(u)-y_3(u))\,p(z,u)|\,du\;dz\, 
\left| \frac{x_2}{{\hat x}{\hat y}}\right|\\
& \phantom{\leq } + \beta  \int_{\zm}^{\zM} \int_{\zm}^{z} |y_3(u)\,p(z,u)|\,du\,dz\, \left|\frac{x_2-y_2}{{\hat y}^2}\right|\\[0.5em]
& \leq \frac{3\beta}{|\hat y|} \, \normx{x-y}.
\end{align*}
Similar relations hold for the other terms $F_1(x,y;w)$, $F_2(x,y;w)$ and $F_3(x,y;w)$. It is then obvious that the norm
$ \|DQ(x;\cdot) -DQ(y;\cdot) \|_{OP}$ tends to zero, for $\normx{x- y}$ going to zero, and the proof is complete.

\bibliographystyle{abbrv}

\end{document}